\input amstex
 \documentstyle{amsppt}
 \input graphicx

\nologo
 \pageno=1
  \loadbold
   \leftheadtext{{\smc Gerd Grubb}} 
    \rightheadtext{{\smc Limited regularity}}

\hsize=5 true in
 \vsize=8.2 true in
  \hoffset=.65 in
   \voffset= 30 pt
    \TagsOnRight
     \NoBlackBoxes

      \font\sevenrm=cmr7


\topinsert\vskip-1.5
 \baselineskip
   {\vbox{\sevenrm\baselineskip 7pt
 \noindent Manuscript submitted to\hfill doi:10.3934/xx.xx.xx.xx \break
  \noindent AIMS journals\hfil\break
   \line {Volume {\sevenbf X}, Number {\sevenbf X},
     XX {\sevenbf 200X}\hfil
     \eightpoint pp. X--XX}}
 }\endinsert

\footnote""{ {\it 2010 Mathematics Subject Classification}.
    Primary: 35K05, 35K35; Secondary 35S11, 47G30, 60G52.}
 \footnote""{ {\it Key words and phrases}.  Fractional Laplacian; stable
process; pseudodifferential operator; fractional heat equation;
fractional Schr\"odinger Dirichlet problem; Lp and H\"older estimates: limited
spatial regularity.}

    \footnote""{   } 

\bigskip
 \document
  \vglue 1\baselineskip

\centerline{\bf  Limited regularity of solutions to fractional heat and Schr\"odinger equations
}

\bigskip
 \medskip

\centerline{\smc Gerd Grubb}
\medskip

{
 \eightpoint
   \centerline{Department of Mathematical Sciences, Copenhagen
 University}
   \centerline{Universitetsparken 5, DK-2100 Copenhagen, Denmark}
 } 
\medskip

\comment
\centerline{\smc First-name2 last-name2}
\medskip

{
 \eightpoint
  \centerline{ First line of the address of the second author}
   \centerline{Other lines, Springfield, MO 65801-2604, USA}
 } 
\endcomment

\bigskip
\centerline{(Communicated by Aim Sciences)}
 \medskip
\bigskip

{
 \eightpoint
  {\narrower\smallskip
   \noindent
    {\bf ABSTRACT.}
When $P$ is the fractional Laplacian $(-\Delta )^a$,
$0<a<1$, or a pseudodifferential generalization thereof, the Dirichlet
problem for the associated heat equation over a smooth set $\Omega \subset{\Bbb R}^n$: $r^+Pu(x,t)+\partial_tu(x,t)=f(x,t)$ on
$\Omega \times \,]0,T[\,$, $u(x,t)=0$ for $x\notin\Omega $, $u(x,0)=0$, is
known to be solvable in relatively low-order Sobolev or H\"older
spaces. We now show that in contrast with differential operator cases,
the regularity of $u$ in $x$ at $\partial\Omega $  when $f$ is very
smooth cannot in general be improved 
beyond a certain estimate. An improvement requires
the vanishing of a Neumann boundary value. --- There is a similar result for
the Schr\"odinger Dirichlet problem $r^+Pv(x)+Vv(x)=g(x)$ on $\Omega $, $\operatorname{supp}
v\subset\overline\Omega  $, with $V(x)\in C^\infty $. The proofs involve a
precise description, of interest in itself, of the
Dirichlet domains in terms of regular functions and
functions with a $\operatorname{dist}(x,\partial\Omega )^a$
singularity.

 \smallskip
  }
   }
    \bigskip

\noindent

\comment
With the setting of the template, automatically the text is set in
10 point fonts, while the abstract and references are in 8 point
fonts. All formulas and pictures must be within the limit of {\bf 5}
inches in width. An abstract is needed and not exceeding {\bf 200}
words. Also needed are some {\bf key words} and {\bf AMS subject
classifications}. Here are some important instructions on how to
prepare your final \TeX\ files. Please pay special attention to the
following:
\endcomment

\head 1. Introduction \endhead


\subhead {1.1 The heat problem} \endsubhead
The main purpose of the paper is to investigate limitations on the regularity of solutions to nonlocal
parabolic Dirichlet problems for $x$ in a bounded smooth subset
$\Omega $ of ${\Bbb R}^n$ and $t$ in an interval $I=\,]0,T[\,$:
$$
\aligned
r^+Pu(x,t)+\partial_tu(x,t)&=f(x,t)\text{ on }\Omega \times I,\\
u(x,t)&=0 \text{ for }x\notin\overline\Omega  ,\\
u(x,0)&=u_0(x),
\endaligned
\tag1.1
$$
where $P$ is the fractional Laplacian $(-\Delta )^a$ on ${\Bbb R}^n$, $0<a<1$, or a
pseudodifferential generalization (that can be $x$-dependent and
non-symmetric); here $r^+$ denotes restriction to $\Omega $. For the
{\it stationary} Dirichlet problem 
$$
r^+Pv(x)=g(x) \text{ in }\Omega ,\quad \operatorname{supp}u\subset \overline\Omega ,\tag1.2
$$
it is known from works of Ros-Oton, Serra, Grubb \cite{47,49,27,26}, that the solution $v$ bears a singularity at
$\partial\Omega $ like $d(x)^a$, where
$d(x)=\operatorname{dist}(x,\partial\Omega )$ near $\partial\Omega $, but that $v/d^a$ is
steadily more regular, the more regular $g$ is. One has for example:
$$
\align
g\in C^\sigma(\overline\Omega )& \implies v/d^a\in C^{a+\sigma}(\overline\Omega )\text{
for }\sigma>0\text{ with }
\sigma ,a+\sigma\notin{\Bbb N},\tag1.3\\
g\in C^\infty (\overline\Omega )& \iff v/d^a\in C^{\infty }(\overline\Omega ).\tag1.4
\endalign$$
(1.3) is shown in \cite{47,49} for small $\sigma $ (allowing low
regularity of $\partial\Omega $), and in \cite{26} for all $\sigma $, and
(1.4) is shown in  \cite{27}, drawing on early work of H\"ormander \cite{35}.

To  formulate results for 
(1.1), let us
temporarily denote  the domain spaces
for the Dirichlet problem (1.2) with  $g\in \overline H_p^r(\Omega )$,
resp.\ $g\in \overline C^r(\Omega )$ by
$$
\align
 D_{r,H_p}(P)&=\{v\in \dot H_p^a(\overline\Omega )\mid
 r^+Pv\in \overline H_p^r(\Omega )\}\text{ for }r\ge 0;\tag1.5\\
 D_{r,C}(P
)&=\{v\in \dot H_2^a(\overline\Omega )\mid
 r^+Pv \in \overline C^r(\Omega )\}\text{ for }r\in {\Bbb R}_+\setminus{\Bbb N}.\tag1.6
\endalign
$$
These spaces, which are known to identify with spaces
$H_p^{a(2a+r)}(\overline\Omega )$, $C_*^{a(2a+r)}(\overline\Omega )$, respetively, introduced in \cite{27,26}, will be
described in a precise way in Section 3, explaining how the factor $d^a$ enters.
For the notation of spaces with dots and overlines, cf.\ (2.4) and
(2.17) below.

In the {\it nonstationary} case, there have up to now been shown some results for (1.1) in function
spaces
of relatively low order, such as e.g., 
under various hypotheses
on $P$: 
$$
\align
&f\text{ is }C^\sigma \text{ in $x$ and $C^{\frac{\sigma }{2a}}$ in $t$}
\implies u/d^a
\text{ is }C^{a+\sigma }\text{ in $x$ and $C^{\frac{a+\sigma }{2a}}$ in
$t$},\tag1.7\\
&f\in L_p(\Omega \times I)
\iff u\in L_p(I; D_{0,H_p}(P))\cap \overline H^1_p(I,L_p(\Omega ))
 ,\tag1.8\\
&f\in L_2(I;\overline H^{r}(\Omega
))\cap \overline H^k(I;L_2(\Omega ))\text{ with }\partial_t^jf(x,0)=0
\text{ for }j<k \tag1.9\\
&\qquad\qquad \implies u \in L_2(I;
D_{r,H_2}(P))\cap \overline
H^{k+1}(I;L_2(\Omega ))
 ,
\endalign$$
with $ \sigma \in \,]0,a],
a+\sigma \notin{\Bbb N}$, $1<p<\infty$, $k\in{\Bbb N}$, $r\le
\min\{2a, a+\frac12-\varepsilon \}$. 
Here (1.7) is shown by Ros-Oton and Vivas \cite{50} building on
Fernandez-Real and Ros-Oton \cite{19}, and  (1.8), (1.9) are shown in
\cite{31,33}.

It is natural to ask whether the nonstationary results can be lifted to
higher regularities in $x$ like in the stationary case: Will solutions have a
$C^\infty $-property if $f$ is $C^\infty $? or e.g.\ a higher H\"older
regularity, when $f$ belongs to a higher H\"older space? Such  rules
holds for differential operator heat problems, and for interior
regularity \cite{31}, but, perhaps
surprisingly, they {\it  do not} hold up to the boundary in the present nonlocal cases.

A first counterexample to the $C^\infty $-lifting was given in
\cite{33},
derived from a certain irregularity of the eigenfunctions of the
Dirichlet realization of $P$ in selfadjoint, $x$-independent cases. In the present
study we show  how the regularity  of the solution is limited to
$u\in D_{a,C}(P
)$ with respect to $x$
(not $D_{a+\delta ,C}(P
)$ with $\delta
>0$), {\it unless the boundary value of $u/d^a$
vanishes}. 

 The heat equation result is based on an analysis of the solutions of the
resolvent equation for $\lambda \ne 0$,
$$
(r^+P-\lambda )v=g \text{ in }\Omega ,\quad \operatorname{supp}v\subset \overline\Omega ,\tag1.10
$$
and a precise description of the spaces
$D_{r,C}(P
)$. 

We also study  Schr\"odinger Dirichlet problems 
$$
(r^+P+V )v=g \text{ in }\Omega ,\quad \operatorname{supp}v\subset \overline\Omega ,\tag1.11
$$
with  a $C^\infty $-potential $V$, and find a related limitation on the
smoothness of solutions, when $V$ does not vanish on $\partial\Omega $.

\subhead{1.2 On the operators $P$}\endsubhead
The fractional Laplacian $(-\Delta )^a$ on ${\Bbb R}^n$, $0<a<1$, 
 can be described as a pseudodifferential operator ($\psi $do) or as
a singular integral operator:
$$
\align
(-\Delta )^au&=\operatorname{Op}(|\xi |^{2a})u=
\Cal F^{-1}(|\xi |^{2a}\hat u(\xi ))\tag1.12\\
&=c_{n,a}PV\int_{{\Bbb R}^n}\frac{u(x)-u(x+y)}{|y|^{n+2a}}\,dy.\tag1.13
\endalign
$$
The operators we shall study are the following generalization of
(1.12) to a large class of $\psi
$do's: The classical strongly elliptic $\psi
$do's $P=\operatorname{Op}(p(x,\xi ))$ of order $2a$,
with 
symbol $p(x,\xi )\sim\sum_{j\in{\Bbb
    N}_0}p_j(x,\xi )$ being {\it even}:
$$
p_j(x,-\xi )=(-1)^jp_j(x,\xi ),\text{ all }j.\tag1.14
$$
An example is $P=A(x,D)^a$,
where $A(x,D)$ is a second-order strongly elliptic differential
operator.

The singular integral definition (1.13) can also be generalized, by
replacement of the kernel function $|y|^{-n-2a}$ by other positive functions
$K(y)$ homogeneous of degree $-n-2a$ and {\it even}, i.e.\ $K(-y)=K(y)$, and with possibly less
smoothness, or nonhomogeneous but estimated 
above in terms of $|y|^{-n-2a}$  (see e.g.\ the survey \cite{45}); this gives translation-invariant
symmetric operators. These are
operators defining stable L\'evy processes. The case where $K$ is
homogeneous and
$C^\infty $ for $y\ne 0$ is a special case of our $\psi $do's, with
symbol $p(\xi )=\Cal F K(y)$.  

Whereas the pseudodifferential symbols $p(x,\xi )$ generally have complex
values and $P$ need not be symmetric, the singular integral definition is usually applied in a real
set-up (with $K$ and $u$ real), defining symmetric operators;
then also $p(\xi )=\Cal F^{-1}K$ is real thanks to the condition
$K(-y)=K(y)$.

Our general assumption is:

\definition{Hypothesis 1.1} For some $a>0$, $P=\operatorname{Op}(p(x,\xi ))$ is a
classical 
$\psi $do 
of order $2a$ on ${\Bbb R}^n$, strongly elliptic (i.e.\
$\operatorname{Re}p_0(x,\xi )\ge c|\xi |^{2a}$ for $|\xi |\ge1$, with $c>0$),  
with even symbol $p(x,\xi )$, cf.\ {\rm (1.14)}, and possibly with a
smoothing term $\Cal R$ added (continuous from $\Cal E'({\Bbb R}^n)$ to
$C^\infty ({\Bbb R}^n)$).
\enddefinition

The $L_2$-realization $P_{\operatorname{Dir},2}$ of the Dirichlet
problem (1.2) is easy to define, as the sectorial operator in
$L_2(\Omega )$ defined by a variational construction from the
associated lower bounded
sesquilinear form $Q_0(u,v)$ with $D(Q_0)=\dot H^a(\overline\Omega )$ 
(defined
by extension by continuity from $(Pu,v)_{L_2(\Omega )}$ on $C_0^\infty
(\Omega )$). Its domain is 
$D(P_{\operatorname{Dir},2})=D_{0,H_2}(P)$ (cf.\ (1.5)), and there is
a Fredholm solvability of (1.2) for $g\in L_2(\Omega )$. The
difficulty lies in describing its domain more precisely, as well as
the domains for other range spaces, as in (1.5)--(1.6). This was done
in \cite{27,26} in terms of the $a$-transmission spaces, recalled in
Section 2.3 below and further studied in Section 3.

\subhead{1.3 Preview of results}\endsubhead
 The
domain spaces $D_{r,H_p}(P)$ (1.5) and  $D_{r,C}(P)$ (1.6) were
 characterized in \cite{27,14} by formulas 
$$
\aligned
D_{r,H_p}(P)&= H_p^{a(2a+r)}(\overline\Omega )\equiv\Lambda _+^{(-a)}e^+\overline H^{a+r}(\Omega
),\\
D_{r,C}(P)&= C_*^{a(2a+r)}(\overline\Omega )\equiv\Lambda
_+^{(-a)}e^+\overline C_*^{a+r}(\Omega )
\endaligned
$$
(explained in Section 2 below); here $\Lambda _+^{(-a)}$ is a special
auxiliary pseudodifferential operator, and the $C^s_*$ denote the
H\"older-Zygmund spaces generalizing the usual H\"older spaces $C^s$,
$s\in{\Bbb R}_+\setminus{\Bbb N}$, to all
real $s$.

Our first task is to give a more elementary description of these
spaces. In fact, the space $\Lambda _+^{(-a)}e^+\overline C_*^{a+r}(\Omega
)$ can be described in terms
of a component in $\dot C^{2a+r}(\overline\Omega )$ (the functions with H\"older smoothness
$2a+r$ in ${\Bbb R}^n$ and support in $\overline\Omega $)  and a
component pulled back from the boundary $\partial\Omega $ as 
$d^a$
times a sum of 
Poisson operators applied to weighted boundary values $\gamma _j^av$,
$$
\gamma _j^av=\Gamma (a+1+j)\gamma _j(v/d^a);\tag1.15
$$
here  $\gamma _ju=(\frac{\partial}{\partial \nu
})^ju|_{\partial\Omega }$.
The description is worked out in Section 3 as a further development of
results from \cite{27,26}. 
With a Poisson operator $K_{(0)}$  
specially constructed for the domain $\Omega $, we find for 
example:

\proclaim{Theorem 1.2}
When $a\in \,]0,1[\,$ and $r>0$ with $a+r\in \,]0,1[\,$
and $2a+r\notin {\Bbb N}$, then
$$
v\in D_{r,C}(P)\iff v=d^aK_{(0)}\psi +w \text{ for }\psi \in
C^{a+r}(\partial\Omega ), w\in \dot C^{2a+r}(\overline\Omega );\tag1.16
$$
here $\psi=\gamma _0(v/d^a)$. 
\endproclaim 

In this context, $\gamma _0(v/d^a)$ (or $\gamma
_0^av$) is usually regarded as a Neumann boundary value, and (1.16) shows how
the irregularity in $v$ comes from the Neumann boundary value.
 For higher
 $r$, with $a+r\in \,]M-1,M[\,$ for a positive integer $M$, there is a
 similar formula with $M$ boundary values.

Similarly, there is a decomposition like (1.16) of $D_{r,H_p}(P)$ in
terms of Sobolev-type spaces when $a+r\in \mathopen{]}\frac1p,1+\frac1p\mathclose{[}$,
and decompositions with more boundary terms for higher $r$.
See details in Theorems 3.4, 3.6, 3.10, 3.12 and their corollaries below. 

Next, by use of such descriptions,
 we show in Sections 4 for the
Schr\"odinger problem:

\proclaim{Theorem 1.3} Let $0<a<1$, and let  $V\in  \overline C^\infty
(\Omega )$. When  $g\in \overline C^{a }(\Omega
)$, the solutions of {\rm (1.11)} are in $D_{a,C}(P)$. 

 Let $V\ne 0$ on an open subset $\Sigma $ of
$\partial\Omega $, and let $v$ be a solution of {\rm (1.11)}. If there is a $\delta >0$ such that  $v\in D_{a+\delta
,C}(P)$ with  $g\in \overline C^{a +\delta }(\Omega
)$,  then $\gamma _0^av$ vanishes on $\Sigma $.

Moreover, if for some noninteger $s>a$, $v\in D_{s
,C}(P)$ with  $g\in \overline C^{s }(\Omega
)$,  then $\gamma _j^av$ vanishes on $\Sigma $ for all integers  $j<s-a$.
\endproclaim

Finally, Section 5 gives the result for the heat problem:

\proclaim{Theorem 1.4} Let $0<a<1$. Let $u(x,t)\in  \overline
W^{1,1}(I;D_{a,C}(P))$ with $u(x,0)=0$; then it satisfies {\rm (1.1)}
with an $f(x,t)\in  L_1(I;\overline C^{a }(\Omega))$. Here $\gamma _0^au\in
W^{1,1}(I;C_*^{2a}(\partial\Omega )) $ and can take any value at a $t>0$.

If there is a $\delta >0$ such that $u(x,t)\in  \overline
 W^{1,1}(I;D_{a+\delta ,C}(P))$  and  $f(x,t)\in $ \linebreak $ L_1(I;\overline C^{a+\delta  }(\Omega))$, 
 then $\gamma _0^au(x,t)= 0$ on $\partial\Omega \times I$.

Moreover, if for some noninteger $s>a$,  $u(x,t)\in  \overline
 W^{1,1}(I;D_{s ,C}(P))$  and  $f(x,t)\in  L_1(I;\overline C^{s }(\Omega))$, 
 then $\gamma _j^au=0$ for all integers  $j<s-a$.
\endproclaim

Expressed in words, a necessary condition for lifting the regularity
parameter
$a$ to $a+\delta $ for $u$ and $f$ is that the Neumann boundary
value vanishes. 

Higher liftings require the vanishing of more traces.

\subhead{1.4 Some references}\endsubhead
There is a rich literature on the fractional Laplacian and its
generalizations and applications in probability, finance, differential
geometry and mathematical physics;
let us mention some of the studies through the times:  
Blumenthal and Getoor \cite{6}, Landkof \cite{42}, Hoh and Jacob
\cite{34},  Kulczycki
\cite{41}, Chen and Song \cite{12}, Jacob \cite{37},
Jakubowski \cite{38}, Bogdan, Burdzy and Chen \cite{7}, Cont and
Tankov \cite{13}, 
Silvestre \cite{52},
Caffarelli and Silvestre \cite{10},  Gonzalez, Mazzeo and
Sire \cite{22}, Musina and Nazarov \cite{44}, Frank and Geisinger
\cite{20}, 
Ros-Oton and Serra \cite{47,48,49}, Abatangelo \cite{1}, Felsinger, Kassmann and Voigt
\cite{18}, Bonforte, Sire
and Vazquez \cite{8}, Servadei and Valdinoci \cite{51}, Ros-Oton
\cite{45,46}, Abatangelo, Jarohs and Saldana \cite{3}. There are many
more papers referred to in these works, and numerous applications to
nonlinear problems.

For the fractional heat equation there are, besides the already
mentioned works \cite{19}, \cite{50}, \cite{31}, \cite{33}, 
contributions to the regularity theory in e.g.\ Felsinger and Kassmann \cite{17}, Chang-Lara and Davila
\cite{11}, Jin and Xiong \cite{39},   Leonori, Peral, Primo and
Soria \cite{43}, 
Biccari, Warma and Zuazua \cite{5}.
For the fractional Schr\"odinger equation, see
e.g.\ Fall \cite{16},  Diaz,
Gomez-Castro and Vazquez \cite{15} and their references.

\head 2.  Preliminaries \endhead

\subhead 2.1 $L_p$ Sobolev spaces \endsubhead 

The domain spaces for homogeneous Dirichlet problems were described
in many scales of spaces in \cite{27} and \cite{26}; we shall here
just focus on two scales, namely the Bessel-potential scale $H^s_p$
($1<p<\infty $) (the main subject of \cite{27}), which serves to show
general estimates in $L_p$ Sobolev spaces, and the H\"older-Zygmund scale 
$C^s_*$ (included in \cite{26} in a systematic way) leading to optimal
H\"older estimates. We shall go through the various concepts somewhat
rapidly, since they have  already been explained in previous
papers; the reader may consult e.g.\ \cite{27,26} if more details
are needed.

We recall that the standard Sobolev spaces $W^{s,p}({\Bbb R}^n)$, $1<p<\infty $ and
$s\ge 0$, have a different character according to whether $s$ is
integer or not. Namely, for $s$ integer, they consist of
$L_p$-functions with derivatives in $L_p$ up to order $s$, hence
coincide with the Bessel-potential spaces $H^s_p({\Bbb R}^n)$, defined
for $s\in{\Bbb R}$ by 
$$
H_p^s({\Bbb R}^n)=\{u\in \Cal S'({\Bbb R}^n)\mid \Cal F^{-1}(\langle{\xi }\rangle^s\hat u)\in
L_p({\Bbb R}^n)\}.\tag2.1
$$
Here $\Cal F$ is the Fourier transform  $\hat
u(\xi )=\Cal F
u(\xi )= \int_{{\Bbb R}^n}e^{-ix\cdot \xi }u(x)\, dx$, and the function $\langle\xi \rangle$ equals $(|\xi |^2+1)^{\frac12}$.
For noninteger $s$, the $W^{s,p}$-spaces coincide with the Besov spaces, defined e.g.\ 
as follows: For $0<s<2$,
$$
f\in
B^s_p({\Bbb R}^{n})\iff \|f\|_{L_p}^p+ \int_{{\Bbb R}^{2n}}\frac{|f(x)+f(y)-2f((x+y)/2)|^p}{|x+y|^{n+ps}}\,dxdy<\infty ;\tag2.2
$$
and $B^{s+t}_p({\Bbb R}^n)=(1-\Delta )^{-t/2}B^s_p({\Bbb R}^n)$ for all $t\in{\Bbb R}$.
The Bessel-potential spaces are important because they are most
directly related to $L_p({\Bbb R}^n)$; the Besov spaces have other
convenient properties, and are
needed for boundary value problems in an $H^s_p$-context,
 because they are the correct range spaces for
trace maps $\gamma _ju=(\partial_n^ju)|_{x_n=0}$:
$$
\gamma _j\colon \overline H^s_p({\Bbb R}^n_+), \overline B^s_p({\Bbb R}^n_+) \to
B_p^{s-j-1/p}({\Bbb R}^{n-1}), \text{ for }s-j-\tfrac1p >0,\tag2.3 
$$
surjectively and with a continuous right inverse; see e.g.\ the overview in
the introduction to \cite{23}. For $p=2$, the two scales are
identical, but for $p\ne 2$ they are related by strict inclusions: $
H^s_p\subset B^s_p\text{ when }p>2$, $H^s_p\supset B^s_p\text{ when
}p<2$. When $p=2$, the index $p$ is usually omitted.

The following subsets of
${\Bbb R}^n$ will be considered:  
 ${\Bbb R}^n_\pm=\{x\in
{\Bbb R}^n\mid x_n\gtrless 0\}$ (where $(x_1,\dots, x_{n-1})=x'$), and
 bounded $C^\infty $-subsets $\Omega $ with  boundary $\partial\Omega $, and
their complements.
Restriction from ${\Bbb R}^n$ to ${\Bbb R}^n_\pm$ (or from
${\Bbb R}^n$ to $\Omega $ resp.\ $\complement\overline\Omega $) is denoted $r^\pm$,
 extension by zero from ${\Bbb R}^n_\pm$ to ${\Bbb R}^n$ (or from $\Omega $ resp.\
 $\complement\overline\Omega $ to ${\Bbb R}^n$) is denoted $e^\pm$. Restriction
 from $\overline{\Bbb R}^n_+$ or $\overline\Omega $ to $\partial{\Bbb R}^n_+$ resp.\ $\partial\Omega $
 is denoted $\gamma _0$. 

We denote by $d(x)$ a function of the form $
d(x)=\operatorname{dist}(x,\partial\Omega )$ for $x\in\Omega $, $x$ near $\partial\Omega $,
extended to a smooth positive function on $\Omega $; $d(x)=x_n$ in the
case of ${\Bbb R}^n_+$.

Along with the spaces $H^s_p({\Bbb R}^n)$ defined in (2.1), we have
the two scales of spaces associated with $\Omega $ for $s\in{\Bbb R}$:
$$
\aligned
\overline H_p^{s}(\Omega)&=\{u\in \Cal D'(\Omega )\mid u=r^+U \text{ for some }U\in
H_p^{s}({\Bbb R}^n)\}, \text{ the {\it restricted} space},\\
\dot H_p^{s}(\overline\Omega )&=\{u\in H_p^{s}({\Bbb R}^n)\mid \operatorname{supp} u\subset
\overline\Omega  \},\text{ the {\it supported} space;}
\endaligned \tag2.4
$$
here $\operatorname{supp}u$ denotes the support of $u$. The definition
is also used with $\Omega ={\Bbb R}^n_+$.  $\overline
H_p^s(\Omega )$ is in other texts often denoted  $H_p^s(\Omega )$  or
$H_p^s(\overline\Omega  )$, and $\dot H_p^{s}(\overline\Omega )$ may be indicated with a
ring, zero or twiddle;
the current notation stems from H\"ormander \cite{36}, Appendix B2.
There are similar spaces with $B^s_p$.

We recall that $\overline H_p^s(\Omega )$ and $\dot H_{p'}^{-s}(\overline\Omega )$ are dual
spaces with respect to a sesquilinear duality extending the $L_2(\Omega )$-scalar
product; $\frac1p+\frac1{p'}=1$.

\subhead 2.2 Pseudodifferential operators \endsubhead

A {\it pseudodifferential operator} ($\psi $do) $P$ on ${\Bbb R}^n$ is
defined from a symbol $p(x,\xi )$ on ${\Bbb R}^n\times{\Bbb R}^n$ by 
$$
Pu=\operatorname{Op}(p(x,\xi ))u 
=(2\pi )^{-n}\int e^{ix\cdot\xi
}p(x,\xi )\hat u\, d\xi =\Cal F^{-1}_{\xi \to x}(p(x,\xi )\Cal F u(\xi
)),\tag 2.5
$$  
using the Fourier transform $\Cal F$. 
We refer to
textbooks such as H\"ormander \cite{36}, Taylor \cite{53}, Grubb \cite{25} for the rules of
calculus (in particular the definition by oscillatory integrals in \cite{36}). The symbols
$p$ of order $m\in{\Bbb R}$ we shall use  are assumed to be {\it
classical}:
$p$  
has an asymptotic expansion $p(x,\xi )\sim \sum_{j\in{\Bbb
N}_0}p_j(x,\xi )$ with $p_j$ homogeneous in $\xi $ of degree $m-j$ for
$|\xi |\ge 1$, all $j$, such that $$
\partial_x^\beta \partial_\xi ^\alpha (p(x,\xi )-
{\sum}_{j<J}p_j(x,\xi )) \text{ is }O(\langle\xi \rangle^{m-\alpha -J})\text{ for
all }\alpha ,\beta \in{\Bbb N}_0^n, J\in{\Bbb N}_0.\tag2.6
$$
 $P$ (and $p$) is  said to be   {\it strongly elliptic}
when $\operatorname{Re}p_0(x,\xi )\ge c|\xi |^m $ for $|\xi |\ge 1$,
with $c>0$. These classical $\psi $do's of order $m$ are continuous from $H^s_p({\Bbb
R}^n)$ to $H^{s-m}_p({\Bbb R}^n)$ for all $s\in{\Bbb R}$. For a
complete theory one adds to these operators 
the smoothing operators (mapping any  $H^s_p({\Bbb R}^n)$ into
$\bigcap_tH^t_p ({\Bbb R}^n)$), regarded as operators of order
$-\infty $. (For example, $(-\Delta )^a$ fits into the calculus when
it is 
written as $\operatorname{Op}((1-\zeta (\xi ))|\xi |^{2a})+\operatorname{Op}(\zeta (\xi
)|\xi |^{2a})$, where $\zeta (\xi )$ is a $C^\infty $-function that
equals $1$ for $|\xi |\le \frac12$ and 0 for $|\xi |\ge 1$; the 
second term is smoothing.)

\example{Remark 2.1}
The operators we consider in this paper are moreover assumed to be {\it even},
cf.\ (1.14), for simplicity. 
In comparison, $P$ of order $2a$ satisfies the $a$-transmission
condition introduced  by H\"ormander
\cite{35,36,27} relative to a given smooth set $\Omega \subset {\Bbb R}^n$ when
$$
       \partial_x^\beta \partial_\xi ^\alpha p_j(x,-\nu(x) ) =(-1)^{|\alpha
       |+j}\partial_x^\beta \partial_\xi ^\alpha p_j(x,\nu (x) ),\text{ all }\alpha ,\beta ,j,\tag2.7
$$
at all points $x\in\partial\Omega $, with interior normal denoted $\nu
(x)$. The evenness means that this is satisfied for {\it any} choice
of $\Omega $. The results in the following hold also when one only
assumes that the
$a$-transmission condition is satified relative to the particular
domain $\Omega $ considered.
\endexample

For our description of the solution spaces for (1.2) we must introduce
 {\it order-reducing
operators}. There is a simple definition of operators $\Xi _\pm^t $ on
${\Bbb R}^n$, $t\in{\Bbb R}$,
$$ 
\Xi _\pm^t =\operatorname{OP}(\chi _\pm^t),\quad \chi _\pm^t(\xi )=(\langle{\xi '}\rangle\pm i\xi _n)^t ;\tag 2.8
$$
 they preserve support
in $\overline{\Bbb R}^n_+m$, respectively, because the symbols extend as holomorphic
functions of $\xi _n$ into ${\Bbb C}_\mp$, respectively; ${\Bbb
C}_\pm=\{z\in{\Bbb C}\mid \operatorname{Im}z\gtrless 0\}$. (The functions
$(\langle{\xi '}\rangle\pm i\xi _n)^t $  satisfy only part of the estimates
(2.6) with $m=t$, but the $\psi $do definition can be applied anyway.)
 There is a more refined choice $\Lambda _\pm^t $
\cite{23, 27}, with
symbols $\lambda _\pm^t (\xi )$ that do
satisfy all the required estimates, and where $\overline{\lambda _+^t }=\lambda _-^{t }$.
These symbols likewise have holomorphic extensions in $\xi _n$ to the complex
halfspaces ${\Bbb C}_{\mp}$, so that
the operators preserve
support in $\overline{\Bbb R}^n_\pm$, respectively. Operators with that property are
called "plus" resp.\ "minus" operators. There is also a pseudodifferential definition $\Lambda
_\pm^{(t )}$ adapted to the situation of a smooth domain $\Omega
$, by \cite{23,27}.

It is elementary to see by the definition of the spaces $H_p^s({\Bbb R}^n)$
in terms of Fourier transformation, that the operators define homeomorphisms 
for all $s$: $
\Xi^t _\pm\colon H_p^s({\Bbb R}^n) \overset\sim\to\rightarrow H_p^{s- t
}({\Bbb R}^n)$.
The special
interest is that the "plus"/"minus" operators also 
 define
homeomorphisms related to $\overline{\Bbb R}^n_+$ and $\overline\Omega $, for all $s\in{\Bbb R}$: 
$$
\aligned
\Xi ^{t }_+\colon \dot H_p^s(\overline{\Bbb R}^n_+ )&\overset\sim\to\rightarrow
\dot H_p^{s- t }(\overline{\Bbb R}^n_+),\quad
r^+\Xi ^{t }_{-}e^+\colon \overline H_p^s({\Bbb R}^n_+ )\overset\sim\to\rightarrow
\overline H_p^{s- t } ({\Bbb R}^n_+ ),\\
\Lambda  ^{(t) }_+\colon \dot H_p^s(\overline\Omega  )&\overset\sim\to\rightarrow
\dot H_p^{s- t }(\overline\Omega ),\quad
r^+\Lambda  ^{(t) }_{-}e^+\colon \overline H_p^s(\Omega  )\overset\sim\to\rightarrow
\overline H_p^{s- t } (\Omega  ),
\endaligned
$$
with similar rules for $\Lambda ^t_\pm$.
Moreover, the operators $\Xi ^t _{+}$ and $r^+\Xi ^{t }_{-}e^+$ identify with each other's adjoints
over $\overline{\Bbb R}^n_+$, because of the support preserving properties.
There is a
similar statement for $\Lambda ^t_+$ and  $r^+\Lambda ^t_-e^+$, and for  $\Lambda ^{(t )}_+$ and $r^+\Lambda ^{(
t )}_{-}e^+$ relative to the set $\Omega $.

In the definition of these operators one can replace $\langle{\xi '}\rangle$ by
$[\xi ']$ (positive, smooth and equal to $|\xi '|$ for $|\xi '|\ge
1$); this is practical when one needs to speak of homogeneous symbols.

\subhead 2.3 The ${\mu}$-transmision spaces \endsubhead

The special {\it ${\mu} $-transmission spaces} were 
introduced for all $\mu \in{\Bbb C}$ by
H\"ormander \cite{35} for $p=2$, cf.\ the account in \cite{27} where
the spaces are redefined and extended to general $p\in \,]1,\infty
[\,$. In the present paper we restrict the attention to
real $\mu >-1$. The spaces are:
$$
\aligned
H_p^{{\mu} (s)}(\overline{\Bbb R}^n_+)&=\Xi _+^{-{\mu} }e^+\overline H_p^{s- {\mu}
}({\Bbb R}^n_+)=\Lambda  _+^{-{\mu} }e^+\overline H_p^{s- {\mu}
}({\Bbb R}^n_+)
,\text{ for }  s> {\mu} -\tfrac1{p'},\\
H_p^{{\mu} (s)}(\overline\Omega )&=\Lambda  _+^{(-{\mu} )}e^+\overline H_p^{s- {\mu}
}(\Omega ),\text{ for }  s> {\mu} -\tfrac1{p'}.
\endaligned\tag 2.9
$$
Observe in particular that
$$
\aligned
H_p^{{\mu}(s)}(\overline\Omega )&=\dot
H_p^s(\overline\Omega ) \text{ for }s-{\mu}\in \,]-\tfrac1{p'},\tfrac1p [\, ,\\
\dot
H_p^s(\overline\Omega )&\subset H_p^{{\mu}(s)}(\overline\Omega )\subset H_{p,\operatorname{loc}}^{s}(\Omega) \text{ for all }s>{\mu}-\tfrac1{p'},
\endaligned
\tag2.10$$
since $e^+\overline H_p^{s- {\mu}}(\Omega )\supset \dot H_p^{s-\mu }(\overline\Omega ) $ for all
$s>{\mu}-\tfrac1{p'}$, with equality if $s-{\mu}\in
\,]-\tfrac1{p'},\tfrac1p [\,$, and since $\Lambda _+^{(-\mu )}$ is elliptic.
Moreover, for $s\ge {\mu}$,
$$
H_p^{{\mu}(s)}(\overline\Omega )\subset H_p^{{\mu}({\mu})}(\overline\Omega )=\dot H_p^{\mu}(\overline\Omega ).\tag2.11 
$$
Recall also from \cite{27} Sect.\ 5 that there is a hierarchy:
$H_p^{{\mu}(s)}(\overline\Omega )\supset H_p^{({\mu}+1)(s)}(\overline\Omega )\supset\dots
\supset  H_p^{({\mu}+j)(s)}(\overline\Omega )$ for $s>{\mu}+j-\tfrac1{p'}$, 
$$
u\in H_p^{({\mu}+j)(s)}(\overline\Omega ) \iff u\in H_p^{{\mu}(s)}(\overline\Omega )\text{ with }\gamma ^{\mu}_0u=\dots=\gamma ^{\mu}_{j-1}u=0,\tag2.12
$$
where (as in (1.15)) $\gamma _j^\mu u=\Gamma (\mu
+1+j)\gamma _j(u/d^\mu )$.
Moreover,
$$
 {\bigcap}_s H_p^{{\mu}(s)}(\overline\Omega )=\Cal E_{\mu}(\overline\Omega )\equiv e^+d^{\mu}\overline C^\infty (\Omega );\tag2.13
$$
the latter space is dense in $H_p^{\mu (s)}(\overline\Omega )$.
It was shown in \cite{27} that
$$
 H_p^{{\mu}(s)}(\overline\Omega )\subset \dot H_p^{s}(\overline\Omega )+e^+d^\mu \overline
 H_p^{s-\mu }(\Omega ),\text{ for }s>\mu  , s-\mu \notin{\Bbb N}.\tag2.14
$$
The inclusion holds with $\dot H_p^{s}(\overline\Omega )$ replaced by $\dot
H_p^{s-\varepsilon }(\overline\Omega )$ if $s-\mu \in{\Bbb N}$.

The great interest of the spaces $H_p^{{\mu}(s)}(\overline\Omega )$ is that they give an exact representation of
the solution spaces for the Dirichlet problem (1.2), and are
independent of $P$. The following
result comes from \cite{27}:

\proclaim{ Theorem 2.2} Let $P$ satisfy Hypothesis
{\rm 1.1}. Consider $v\in \dot H^{a-1/p'+\varepsilon }_p(\overline\Omega )$ (any $\varepsilon >0$). For $r>-a-\tfrac1{p'}$, the
solutions of problem {\rm (1.2)} with
$g\in \overline H^r_p(\Omega )$ satisfy  $v\in H_p^{a(2a+r)}(\overline\Omega )$; in fact
$$
g\in \overline H^r_p(\Omega )\iff v\in H_p^{a(2a+r)}(\overline\Omega ),\tag2.15
$$
 and the mapping $r^+P\colon v\mapsto g$
is Fredholm between these spaces.

It follows that the Dirichlet domain  $D_{r,H_p}(P)$ defined in {\rm
(1.5)} for $r\ge 0$ satisfies
$$
D_{r,H_p}(P)= H_p^{a(2a+r)}(\overline\Omega ).\tag2.16
$$

\endproclaim

\demo{Proof} The first statement is a case of Th.\ 4.4 of
\cite{27}, with $\mu _0=a$, $m=2a$, the factorization index being $a$
because of the strong ellipticity (as in Eskin [E71]; a detailed
discussion of factorizations is given in \cite{30}).

The second statement specializes this to $r\ge 0$ (which is all we
need in the present paper) where the prerequisite on $v$ can be
simplified to $v\in \dot H_p^a(\overline\Omega )$, since $H_p^{a(2a+r)}(\overline\Omega 
)\subset \dot H_p^a(\overline\Omega )\subset \dot H_p^{a-1/p'+\varepsilon
}(\overline\Omega )$ for small $\varepsilon $. \qed 
\enddemo

\subhead 2.4 Generalizations to H\"older-Zygmund spaces \endsubhead

In \cite{26} (that was written after \cite{27}), the results were extended to many other scales of spaces,
such as Besov spaces $B^s_{p,q}$ for $1\le p,q\le \infty $,  and
Triebel-Lizorkin spaces $F^s_{p,q}$ for $1\le p< \infty $, $1\le q\le \infty $. Of
particular interest is the scale $B^s_{\infty ,\infty }$, also denoted
$C^s_*$, the H\"older-Zygmund scale. Here $C^s_*$ identifies with the
H\"older space $C^s$ when $s\in
{\Bbb R}_+\setminus {\Bbb N}$, and for positive integer $k$ satisfies
$ C^{k-\varepsilon }\supset C^k_*\supset 
C^{k-1,1}\supset  C^k$ for small $\varepsilon >0$; moreover,
$C^0_*\supset L_\infty \supset C^0$. Similarly to (2.4) we denote the
spaces of restricted, resp.\ supported distributions
$$
\aligned
\overline C_*^{s}(\Omega)&=\{u\in \Cal D'(\Omega )\mid u=r^+U \text{ for some }U\in
C_*^{s}({\Bbb R}^n)\},\\
\dot C_*^{s}(\overline\Omega )&=\{u\in C_*^{s}({\Bbb R}^n)\mid \operatorname{supp} u\subset
\overline\Omega  \},
\endaligned
\tag2.17
$$
where
 the star can be omitted when $s\in {\Bbb R}_+\setminus {\Bbb N}$. The $\mu
 $-transmission spaces are defined as 
$$
C_*^{\mu  (s)}(\overline\Omega )=\Lambda  _+^{(-{\mu} )}e^+\overline C_*^{s- {\mu}
}(\Omega ),\text{ for }  s> {\mu} -1,
\tag2.18
$$
and
there are inclusions as described for $H^{{\mu}(s)}(\overline\Omega )$-spaces in
(2.10)--(2.14). In particular:
$$
\aligned
C_*^{{\mu}(s)}(\overline\Omega )&=\dot
C_*^s(\overline\Omega ) \text{ for }s-{\mu}\in \,]-1,0 [\, ,\\
\dot
C_*^s(\overline\Omega )& \subset C_*^{{\mu}(s)}(\overline\Omega )\subset C_{*,\operatorname{loc}}^s(\Omega )\text{ for all }s>{\mu}-1,\\
u\in C_*^{({\mu}+j)(s)}(\overline\Omega )& \iff u\in
C_*^{{\mu}(s)}(\overline\Omega )\text{ with }
\gamma ^{\mu}_0u=\dots=\gamma
^{\mu}_{j-1}u=0,\\
 {\bigcap}_s C_*^{{\mu}(s)}(\overline\Omega )&=\Cal E_{\mu}(\overline\Omega ),\\
C_*^{{\mu}(s)}(\overline\Omega )&\subset \dot C_*^{s}(\overline\Omega )+e^+d^\mu \overline
 C_*^{s-\mu }(\Omega ),\text{ for }s>\mu  , s-\mu \notin{\Bbb N};
\endaligned\tag2.19
$$ 
in the third line $s>\mu +j-1$ is assumed, and
the last line holds with $\dot C_*^{s}(\overline\Omega )$ replaced by $\dot
C_*^{s-\varepsilon }(\overline\Omega )$ if $s-\mu \in{\Bbb N}$.
There is a result similar to that of Theorem 2.2 with $H^s_p$-spaces replaced
by $C^s_*$-spaces:

\proclaim{ Theorem 2.3} Let $P$ satisfy Hypothesis
{\rm 1.1}. 
Consider $v\in \dot H^{a-1/p'+\varepsilon }(\overline\Omega )$ (any
$\varepsilon >0$). For $r\ge 0$, the
solutions of problem {\rm (1.2)} with
$g\in \overline C_*^r(\Omega )$ satisfy  $v\in C_*^{a(2a+r)}(\overline\Omega )$; in fact
$$
g\in \overline C_*^r(\Omega )\iff v\in C_*^{a(2a+r)}(\overline\Omega ),\tag2.20
$$
 and the mapping $r^+P\colon v\mapsto g$
is Fredholm between these spaces.

It follows that the Dirichlet domain  $D_{r,C}(P)$ defined as
$
 D_{r,C}(P)=\{v\in \dot H_2^a(\overline\Omega )\mid
 r^+Pv \in \overline C^r_*(\Omega )\}$ for $r\ge 0$, satisfies
$$
D_{r,C}(P)= C_*^{a(2a+r)}(\overline\Omega ).\tag2.21
$$
\endproclaim

For $r\in{\Bbb R}_+\setminus {\Bbb N}$, this is the domain defined in
(1.6); and when $a+r $ and $2a+r$ are noninteger, $
C_*^{a(2a+r)}(\overline\Omega )$ identifies with $ C^{a(2a+r)}(\overline\Omega )=\Lambda  _+^{(-a )}e^+\overline C^{a+r
}(\Omega )$ defined in terms of ordinary H\"older spaces. 
It is sometimes an
advantage to keep the $C_*$-notation, since one does not have to make
exceptions for integer indices all the time.

Note also that since $ \overline C_*^r(\Omega )\supset  \overline C^r(\Omega )$ when
$r\in{\Bbb N}_0$, we have a regularity implication $g\in \overline
C^r(\Omega )\implies v\in C_*^{a(2a+r)}(\overline\Omega )$ when $r\in{\Bbb
N}_0$.

In applications of the above results it is important to get  a better
understanding of what the ${\mu}$-transmission spaces consist of. Such an
analysis was carried out in local coordinates for the scale $ H_p^{\mu (t)}$ in
\cite{27} and for $ C_*^{\mu (t)}$ in \cite{26}, and in the
next section we take it up again, showing very explicit global results relative
to the set $\Omega $.

\head 3. Analysis of $\mu$-transmission spaces \endhead

\subhead 3.1 Decompositions in terms of the first trace
\endsubhead

In the following, we shall give a characterization of the
$\mu$-transmission spaces  
showing an exact decomposition of the elements in the case of a general
domain $\Omega $; it is a further development of the
decomposition principle described in Th.\ 5.4 of \cite{27}. For clarity, we begin with a
decomposition with just one trace involved.
Recall that $\mu $ is real $>-1$.

In the proofs we shall use a localization with particularly convenient
coordinates, described in detail in \cite{29} Rem.\ 4.3 and Lem.\ 4.4
and recalled in \cite{32} Rem.\ 4.3, which we also recall here:

\example{Remark 3.1}
$\overline\Omega $ has a finite cover by bounded open sets $U_0,\dots, U_{I}$ with
$C^\infty $-dif\-fe\-o\-mor\-phisms
$\kappa _i\colon U_i\to V_i$, $V_i$ bounded open in ${\Bbb R}^n$, such
 that $U_i^+=U_i\cap \Omega $ is mapped to $V_i^+=V_i\cap{\Bbb R}^n_+$ and
 $U_i'=U_i\cap\partial\Omega $ is mapped to $V_i'=V_i\cap\partial\overline{\Bbb R}^n_+$; as
 usual we write $\partial\overline{\Bbb R}^n_+={\Bbb R}^{n-1}$. 
For any such cover there exists an associated partition of unity, namely a family
of functions $\varrho _i\in C_0^\infty (U_i)$ taking values in $[0,1]$
such that $\sum_{i=0,\dots,I}\varrho _i$ is 1 on a neighborhood of $\overline\Omega $.

When $P$ is a $\psi
 $do on ${\Bbb R}^n$, its application to functions
 supported in $U_i$ carries over to functions on $V_i$ as a $\psi $do
 $\underline P^{(i)}$ defined by 
$$
\underline P^{(i)}v=P(v\circ \kappa _i)\circ \kappa _i^{-1},\quad v\in C_0^\infty (V_i).\tag3.1
$$
When $u$ is a function on $U_i$, we usually denote the resulting function $u\circ \kappa
_i^{-1}$ on $V$ by $\underline u$.

We shall use a convenient system of coordinate charts as described 
in \cite{29}, Remark
4.3: Here $\partial\Omega $ is covered with coordinate charts  $\kappa '_i\colon U'_i\to V'_i\subset {\Bbb R}^{n-1}$,
$i=1,\dots, I$, and the $\kappa _i$ will be defined on certain subsets  of a tubular
neighborhood $\Sigma _r=\{x'+t\nu (x')\mid x'\in\partial\Omega , |t|<r\}$,
where $\nu (x')=(\nu _1(x'),\dots,\nu _n(x'))$ is the interior normal
to $\partial\Omega $ at $x'\in\partial\Omega $, and $r$ is taken so
small that the mapping $ x'+t\nu (x')\mapsto (x',t)$ is a
diffeomorphism from $\Sigma _r$ to $\partial\Omega \times
\,]-r,r[\,$. For each $i$, $\kappa _i$ is defined as the mapping $\kappa
_i\colon x'+t\nu (x')\mapsto (\kappa '_i(x'),t)$ ($x'\in U'_i$).  $\kappa _i$ goes from $U_i$ to
$V_i$, where
$$
U_i=\{x'+t\nu (x')\mid x'\in U'_i, |t|<r\}, \quad V_i=V'_i\times
\,]-r,r[\,.\tag3.2
$$
These charts are supplied with a chart
consisting of the identity mapping on an open set $U_0$ containing $\Omega
\setminus \overline \Sigma _{r,+}$, with $\overline U_0\subset \Omega
$, to get a full cover of $\overline\Omega $.

Note that the normal $\nu (x')$ at $x
'\in\partial\Omega $ is carried over to the normal $(0,1)$ at $(\kappa
'_i(x'),0)$ when $x'\in U_i'$. The halfline $L_{x'}=\{x'+t\nu (x')\mid t\ge 0\}$ is the
geodesic into $\Omega $ orthogonal to $\partial\Omega $ at $x'$ (with
respect to the Euclidean metric on ${\Bbb R}^n$), and there is a positive
$r'\le r$ such that for $0<t<r'$, the distance $d(x)$ between  $x=x'+t\nu
(x')$ and $\partial\Omega $
equals $t$. Then $t$ plays the role of $d$ in the definition of
expansions and boundary values of $u\in \Cal E_{\mu }(\overline\Omega )$
 in
\cite{27} (5.3)ff.: 
$$
u=\tfrac 1{\Gamma (\mu +1)}t^{\mu }u_0+\tfrac 1{\Gamma (\mu +2)}t^{\mu +1}u_1+\tfrac
1{\Gamma (\mu +3)}t^{\mu +2}u_2+\dots\text{ for }t>0,\quad u=0\text{ for }t<0,\tag3.3
$$
where the $u_j$ are constant in $t$ for $t<r'$;  this serves to define
the boundary values 
$$
\gamma ^{\mu }_{j}u=\gamma _0u_j\,(=u_j|_{t=0}),\quad
j=0,1,2,\dots\tag3.4
$$
(denoted $\gamma _{\mu ,j}u$ in \cite{27}).

By comparison of (3.3) with
$t^{\mu }$ times the
Taylor expansion of $u/t^{\mu }$ in $t$, we also have:
$$
\gamma _0^{\mu }u=\Gamma (\mu +1)\gamma _0(u/t^{\mu }),\quad \gamma
_1^{\mu }u=\Gamma (\mu +2)\gamma _1(u/t^{\mu })=\Gamma (\mu +2)\gamma
_0(\partial_t(u/t^{\mu })),\text{ etc.}\tag3.5
$$

\endexample

We first recall (and reprove) a result from \cite{27} for the case
where the domain is ${\Bbb R}^n_+$.

\example{Remark 3.2} There is a notational ambiguity in the fact
that the functions we deal with on $\overline{\Bbb R}^n_+$ are often understood as
extended by zero on ${\Bbb R}^n_-$. Here $\gamma _j$ and $\gamma ^\mu _j$ will
always be read as taking boundary values from the interior
${\Bbb R}^n_+$. This could be underlined by writing $\gamma _j, 
\gamma ^\mu _j$ as $\gamma
_j^+, \gamma _j^{\mu ,+}$, but we refrain from this notational complication. 

The Poisson operators $K$ in the Boutet de Monvel calculus map spaces over
${\Bbb R}^{n-1}$ into spaces over ${\Bbb R}^n_+$, but in their concrete definition
by Fourier transformation $K\varphi =$\linebreak $\Cal F^{-1}_{\xi \to
x}(k(x',\xi )\hat\varphi (\xi '))$, the symbols $k(x',\xi )$ have
$\Cal F^{-1}_{\xi _n\to x_n}k$ supported for $x_n\ge 0$,  whereby 
$K$ actually maps into the function spaces extended by zero on
${\Bbb R}^n_-$. It is customary to leave out the explicit mention of $e^+$ (to
read $K$ as $e^+K$ if the context requires it). The indication $e^+$
can be mentioned to underline the mapping property, but will most often
be left out.

Similar principles are followed when ${\Bbb R}^n_+$ is replaced by $\Omega $.
\endexample

\proclaim{Lemma 3.3} 

Let $K_0$ denote the Poisson operator from ${\Bbb
R}^{n-1}$ to $\overline{\Bbb R}^n_+$ with symbol $(\langle{\xi '}\rangle+i\xi _n)^{-1}$.

When   $s>\mu +\tfrac1p $, the elements of $H_p^{\mu (s)}(\overline{\Bbb R}^n_+)$ have a
unique decomposition
$$
\aligned
u&=v +w ,\text{ where }w \in H_p^{(\mu +1)(s)}(\overline{\Bbb R}^n_+),\text{ and}\\
v &=e^+\tfrac1{\Gamma (\mu +1)}x_n^{\mu } K_0 \gamma ^\mu _0u\in e^+x_n^\mu \overline H_p^{s-\mu }({\Bbb R}^n_+ )\cap H_p^{\mu (s)}(\overline{\Bbb R}^n_+).
\endaligned
\tag 3.6 
$$
In fact, the elements of $H_p^{\mu (s)}(\overline{\Bbb R}^n_+)$ are parametrized as 
$$
u=e^+\tfrac1{\Gamma (\mu +1)}x_n^{\mu } K_0 \varphi   +w ,\tag3.7
$$
where $\varphi  $ runs through $B_p^{s-\mu -1/p}({\Bbb R}^{n-1} )$ and
$w $ runs through $ H_p^{(\mu +1)(s)}(\overline{\Bbb R}^n_+)$; here $\varphi  $ equals $\gamma ^\mu _0u$.

\endproclaim 

\demo{Proof} 
In detail, $K_0$ is
the elementary Poisson operator of order $0$ in the Boutet de Monvel
calculus
(cf.\ e.g.\  \cite{9,24,25}):
$$
K_0\colon \varphi \mapsto \Cal F^{-1}_{\xi \to x}(\hat \varphi (\xi ')
(\langle{\xi '}\rangle +i\xi _n)^{-1})= \Cal F^{-1}_{\xi '\to x'}(\hat \varphi (\xi ')e^+r^+e^{-x_n\langle{\xi '}\rangle }),\tag3.8
$$
which solves the Dirichlet problem $(1-\Delta )u=0$ on ${\Bbb R}^n_+$, $\gamma
_0u=\varphi $ on ${\Bbb R}^{n-1}$. 
(We use conventions as in
Remark 3.2.)
Define $K^\mu _0$  by
$$
\aligned
K^\mu _0\varphi &=\Xi ^{-\mu }_+e^+K_0\varphi = \Cal F^{-1}_{\xi \to x}(\hat
\varphi (\xi ')(\langle{\xi '}\rangle +i\xi _n)^{-1-\mu })\\
&=\tfrac1{\Gamma (\mu +1)}x_n^\mu  \Cal F^{-1}_{\xi '\to x'}(\hat \varphi (\xi ')e^+r^+e^{-x_n\langle{\xi '}\rangle })=\tfrac1{\Gamma (\mu +1)}x_n^\mu  e^+K_0\varphi ;
\endaligned\tag3.9
$$
by the last expression, it is a right inverse of $\gamma _0^\mu \colon
u\mapsto \Gamma (\mu +1)\gamma _0(u/x_n^\mu )$. 
(These calculations played an important role in \cite{27},
 cf.\ (3.5) and the proofs of
Cor.\ 5.3 and Th.\ 5.4 there. The constant called  $c_\mu$ in (5.16)
there is
written explicitly here, equal to $1/ \Gamma (\mu +1)$.)

When $u\in H_p^{\mu (s)}(\overline{\Bbb R}^n_+)$ for some $s>\mu +\tfrac1p $, then $\gamma _0^\mu u\in
B_p^{s-\mu -1/p}({\Bbb R}^{n-1})$ (cf.\ \cite{27} Th.\ 5.1),  and $w =u-K^\mu _0\gamma _0^\mu u$ has $\gamma
_0^\mu w =0$.  By the mapping properties of the Poisson operator
$K_0$ known from \cite{23}, $e^+K_0\gamma _0^\mu u$
lies in $e^+\overline H_p^{s-\mu }({\Bbb R}^n_+)$. Then the last expression for $K^\mu _0 $ in
(3.9) shows that  $K^\mu _0\gamma _0^\mu u\in x_n^\mu e^+\overline H_p^{s-\mu }({\Bbb R}^n_+)$. 
 Moreover, $K^\mu _0\gamma _0^\mu u$ lies in $H_p^{\mu (s)}(\overline{\Bbb R}^n_+)$,
since it is $\Xi _+^{-\mu }$ of a function in $e^+\overline H_p^{s-\mu
}({\Bbb R}^n_+)$ (by the first equality in (3.9) with $\varphi =\gamma _0^\mu u$).
Then also $w $ lies in $H_p^{\mu (s)}(\overline{\Bbb R}^n_+)$, with $\gamma _0^\mu
w =0$. That $w$ is in the subspace $
H_p^{(\mu +1)(s)}(\overline{\Bbb R}^n_+)$ follows from (2.12) recalled above; we can
also argue more directly as follows: Let $u_n\to u$ in $H_p^{\mu (s)}(\overline{\Bbb R}^n_+)$,  $u_n\in \Cal E_\mu (\overline{\Bbb R}^n_+)\cap
\Cal E'({\Bbb R}^n)$; then $w_n=u_n-K^\mu _0\gamma _0^\mu u_n$ converges to $w$ in
$H_p^{\mu (s)}(\overline{\Bbb R}^n_+)$. Here $w_n$ is in $ \Cal E_\mu (\overline{\Bbb R}^n_+)$ and has
$\gamma _0^\mu w_n=0$, hence lies in $\Cal E_{\mu +1}(\overline{\Bbb R}^n_+)$ (compare
Taylor expansions at $x_n=0$); then the
limit $w$ lies in the closed subspace $H_p^{(\mu +1)(s)}(\overline{\Bbb R}^n_+)$ (cf.\ Prop.\ 4.3 of \cite{27}) of
$H_p^{\mu (s)}(\overline{\Bbb R}^n_+)$. 

The decomposition is unique since $w$ is determined from $u$. 

All functions $\varphi  
\in B_p^{s-\mu -1/p}({\Bbb R}^{n-1})$ give rise to functions in $
H_p^{\mu (s)}(\overline{\Bbb R}^n_+)$ by the mapping $\Xi _+^{-\mu }e^+K_0$, and 
$H_p^{(\mu +1)(s)}(\overline{\Bbb R}^n_+)\subset H_p^{\mu (s)}(\overline{\Bbb R}^n_+)$,
so all functions  $\varphi  
\in $ \linebreak $ B_p^{s-\mu -1/p}({\Bbb R}^{n-1})$  and  $w \in
H_p^{(\mu +1)(s)}(\overline{\Bbb R}^n_+)$ are reached in the decomposition. \qed 


\enddemo

We shall now show a similar result for arbitrary smooth bounded sets $\Omega $.

\proclaim{Theorem 3.4} Let $\Omega \subset {\Bbb R}^n$, bounded, open
with $C^\infty $-boundary. 

There is a
Poisson operator $ K_{(0)}$ of order $0$
from $\partial\Omega $ to $\overline\Omega $ 
(in the Boutet de Monvel calculus) with principal symbol $([\xi
']+i\xi _n)^{-1}$ in local coordinates at the boundary,
such that $K_{(0)}$ is a right inverse of $\gamma _0$, and the
following holds:

The operator $K^\mu _{(0)}$ defined by
$$
 K^\mu _{(0)}= \tfrac1{\Gamma (\mu +1)}d^\mu e^+ K_{(0)},\tag3.10
$$
 maps  $B_p^{s-\mu -1/p}(\partial\Omega )$ into $ e^+d^\mu \overline
H_p^{s-\mu }(\Omega )\cap H_p^{\mu (s)}(\overline\Omega )$ for $s>\mu +\tfrac1p $,
and satisfies
$$
\gamma ^\mu _0 K^\mu _{(0)}\varphi =\varphi,\text{ all }\varphi \in B_p^{s-\mu -1/p}(\partial\Omega ) .\tag 3.11
$$
When $s>\mu +\tfrac1p $, the elements of $H_p^{\mu (s)}(\overline\Omega )$ have a
unique decomposition
$$
u=K^\mu _{(0)} \varphi     +w ,\tag 3.12
$$
where $\varphi    $ runs through $B_p^{s-\mu -1/p}(\partial\Omega )$ and
$w  $ runs through $ H_p^{(\mu +1)(s)}(\overline\Omega )$ (equal to $\dot
H_p^s(\overline\Omega )$ when $s-\mu \in\,]\tfrac1p ,1+\tfrac1p [\,$); here  $\gamma ^\mu
_0u$ equals $\varphi    $.

 \endproclaim

\demo{Proof} 
We use the local coordinates $\kappa _i\colon U_i\to V_i$, $i=0,1,\dots, I$, described in
Remark 3.1, with an associated partition of unity  $\{\varrho
_i\}_{i\le I}$.  We can
moreover choose nonnegative functions $\psi ^k _i\in C_0^\infty
(U_i)$, $k=1,2,3$, such that
$\psi ^1 _i\varrho _i=\varrho _i$, i.e., $\psi ^1_i$ is 1 on $\operatorname{supp} \varrho
_i$, and similarly $\psi ^2_i\psi ^1_i=\psi ^1_i$, $\psi ^3_i\psi ^2_i=\psi ^2_i$.

Let $u\in H_p^{\mu (s)}(\overline\Omega )$, i.e., $u=\Lambda _+^{(-\mu )}z$ for some
$z\in e^+\overline H_p^{s-\mu }(\Omega )$. Write
$$
u=\Lambda _+^{(-\mu )}{\sum}_{i=0}^I\varrho _iz=
{\sum}_i\psi ^1_i\Lambda _+^{(-\mu )}\varrho _iz+{\sum}_i(1-\psi ^1_i)\Lambda _+^{(-\mu )}\varrho _iz.\tag3.13
$$
Since $(1-\psi ^1 _i)\varrho _i=0$, $(1-\psi ^1 _i)\Lambda _+^{(-\mu )}\varrho
_i$ is a $\psi $do of order $-\infty $, so it maps $z$ into $C^\infty
({\Bbb R}^n)$; moreover, its symbol in local coordinates is holomorphic
for $\operatorname{Im}\xi _n<0$, so it preserves support in $\overline\Omega $.
Hence the terms in the second sum in the right-hand side of (3.13) are in $\dot C^\infty (\overline\Omega )$,
contained in $\dot H^t_p(\overline\Omega )\subset H_p^{\mu (t)}(\overline\Omega )$ for all
$t$ (and absorbed in the $w$-term in the final formula). Henceforth we can focus on
the first sum
$$
{\sum}_iu_i,\quad u_i=\psi ^1_i\Lambda _+^{(-\mu )}\varrho _iz;
$$
where $u_i$ is compactly supported in the set $U_i$ and belongs to
$H_p^{\mu (s)}(\overline\Omega )$.

Consider one $u_i$. Here $\underline u_i=u_i\circ \kappa _i^{-1}$ is
in $H_p^{\mu (s)}(\overline{\Bbb R}^n_+)$ with support in $\operatorname{supp}\underline
\psi ^1 _i$.  

By Lemm 3.3,
$$
\underline u_i=e^+\tfrac1{\Gamma (\mu +1)}x_n^\mu  K_0\gamma _0^\mu \underline
u_i+\underline w_i, \quad \underline w_i\in  H_p^{(\mu +1)(s)}(\overline\Omega ),
$$
where $\gamma ^\mu _0\tfrac1{\Gamma (\mu +1)}x_n^\mu  K_0=I$.
Multiplication by $\underline \psi ^2 _i$ or $\underline\psi ^3 _i$ does not alter $\underline
u_i$; this gives the representation (where we denote $\gamma _0\psi ^k_i=\psi ^k_{i,0}$)
$$
\underline u_i=e^+\tfrac1{\Gamma (\mu +1)}\underline \psi ^3 _ix_n^\mu 
K_0\gamma _0^\mu (\underline \psi ^2_i\underline
u_i)+\underline \psi ^3 _i\underline w_i=e^+\tfrac1{\Gamma (\mu +1)}x_n^\mu \underline \psi ^3 _i
K_0\underline \psi ^2 _{i,0}\gamma _0^\mu \underline
u_i+\underline \psi ^3 _i\underline w_i.
$$
Here $e^+\underline \psi ^3 _ix_n^\mu 
K_0\gamma _0^\mu (\underline \psi ^2_i\underline
u_i)\in H_p^{\mu (s)}(\overline{\Bbb R}^n_+)$, since it is compactly supported and is the sum of $e^+x_n^\mu 
K_0\gamma _0^\mu (\underline \psi ^2_i\underline
u_i)\in H_p^{\mu (s)}(\overline{\Bbb R}^n_+)$ and $(1-\underline \psi ^3 _i)e^+x_n^\mu 
K_0\gamma _0^\mu (\underline \psi ^2_i\underline
u_i)\in \Cal E_a(\overline{\Bbb R}^n_+)$, using that $(1-\underline\psi ^3
_i)K_0\underline\psi ^2_{i,0}$ is a Poisson operator of order $-\infty $. Then also
$\underline\psi ^3 _i\underline w_i$ is in  $H_p^{\mu (s)}(\overline{\Bbb R}^n_+)$, and
since its first boundary value $\gamma ^\mu _0(\underline\psi ^3
_i\underline w_i)=
\underline\psi ^3_{i,0}\gamma ^\mu _0\underline w_i$ vanishes, it is in fact
in $ H_p^{(\mu +1)(s)}(\overline{\Bbb R}^n_+)$. 

Denote
$(\underline \psi ^3_i\underline w_i)\circ \kappa _i=\tilde w_i$, then
we get the formula in the original coordinates:
$$
u_i=e^+\tfrac1{\Gamma (\mu +1)}d^\mu  \tilde K^i_0 \gamma _0^\mu 
u_i+\tilde w_i,\text{ where }\tilde K^i_0=(\underline \psi ^3 _i
K_0\underline\psi ^2_{i,0})^\sim,
$$
the operator induced by $\underline\psi ^3 _i K_0 \underline\psi ^2
_{i,0}$ in the original coordinates. Similarly as before, 
 $\gamma
^\mu _0\tfrac1{\Gamma (\mu +1)}d^\mu \tilde K^i_0\gamma _0^\mu u_i=\gamma
_0^\mu u_i$. Finally we find by summation over $i$
the formula
$$
\aligned
u&=e^+\tfrac1{\Gamma (\mu +1)}d^\mu  K_{(0)} \gamma _0^\mu 
u+ w= K^\mu _{(0)} \gamma _0^\mu 
u+ w,\\
\text{ with }K_{(0)}&={\sum}_{i=0}^{I}\tilde K^i_0,\;
K^\mu _{(0)}=e^+\tfrac1{\Gamma (\mu +1)}d^\mu  K_{(0)} \text{ and
}w={\sum}_i\tilde w_i;
\endaligned 
$$
here
$$
\gamma _0K_{(0)}=I\text{ and } \gamma ^\mu _0K^\mu _{(0)}=I.
$$
This shows the asserted unique decomposition, and the mapping
properties follow from those of the localized pieces.
\qed

An immediate corollary is a description of the domain space for
the Dirichlet realization of $P$ in $L_p(\Omega )$, when $0<a<1$ (leaving out the normalization by $\Gamma
$-factors):

\proclaim{Corollary 3.5} When $P$ satisfies Hypothesis {\rm 1.1}
with  $0<a<1$, the domain
$$
D(P_{\operatorname{Dir},p})=\{v\in \dot H_p^a(\overline\Omega )\mid
 r^+Pv\in L_p(\Omega )\}\tag3.14
$$
(called $D_{0,H_p}(P)$ in {\rm (1.5)}) satisfies: $D(P_{\operatorname{Dir},p})=\dot H_p^{2a}(\overline\Omega )$ when
$a<\tfrac1p $, and 
$$
D(P_{\operatorname{Dir},p})=\{u=d^aK _{(0)} \psi +w \mid \psi
\in B_p^{a -1/p}(\partial\Omega ), w\in \dot H_p^{2a}(\overline\Omega )\},\tag3.15
$$
when $a\in \,]\tfrac1p ,1[\,$; here $\psi =\gamma _0(u/d^a)$.
\endproclaim

As shown in \cite{26}, the results of \cite{27} carry over to many other interesting
scales of function spaces. We shall here in particular consider the
H\"older-Zygmund scale, for which we get the following version of
Theorem 3.4:

\enddemo

\proclaim{Theorem 3.6} With the Poisson operator  $K_{(0)}$ defined as in Theorem
{\rm 3.4}, and $K^\mu _{(0)}$ defined by {\rm (3.10)}, the following holds:

For $s>\mu $,  $K^{\mu }_{(0)}$ maps $C_*^{s-\mu }(\partial\Omega )\to
e^+d^{\mu }\overline C_*^{s-\mu }(\Omega )\cap C_*^{(\mu )(s)}(\overline\Omega )$.
Moreover, the elements of
 $C^{\mu (s)}_*(\overline\Omega )$ have unique decompositions
$$
u= K^\mu _{(0)}\varphi  +w, \text{ where }
\tag3.16$$
$\{\varphi ,w\}$ runs through $\in C_*^{s-\mu }(\partial\Omega )\times
C_*^{(\mu +1)(s)}(\overline\Omega )$;   here $\varphi $ equals $\gamma ^\mu _0u$.

The space
$C_*^{(\mu +1)(s)}(\overline\Omega )$ equals $\dot C^s_*(\overline\Omega )$ if $s-\mu \in \,]0,1[\,$.
\endproclaim

\demo{Proof} We use that $C^s_*=B^s_{\infty ,\infty }$ in the Besov
scales $B^s_{p,q}$, where the $\psi $do's and the boundary operators from the
Boutet de Monvel calculus act similarly as in $H^s_p$, as shown in
Johnsen \cite{40}, the consequences for our calculations in the Besov scales being
recalled in \cite{26}. Here $p=\infty $, so $p'=1$. The statements
$s>\mu +\tfrac1p $ and $s-\mu \in \,]\tfrac1p ,1+\tfrac1p [\,$ are here replaced by
$s>\mu $ and $s-\mu \in \,]0,1[\,$.
\qed
\enddemo

For example, this gives: 

\proclaim{Corollary 3.7} When $P$ satisfies Hypothesis {\rm 1.1}
with  $0<a<1$, the domain
$D_{r,C}(P)$ defined in  {\rm (1.6)} satisfies, if $r>0$ with $a+r\in
\,]0,1[\,$, $2a+r\notin {\Bbb N}$:
$$
u\in D_{r,C}(P)\iff u=d^aK_{(0)}\psi +w \text{ for }\psi \in
C^{a+r}(\partial\Omega ), w\in \dot C^{2a+r}(\overline\Omega );
$$
here $\psi =\gamma _0(u/d^a)$.
\endproclaim

This shows how the irregularity of $u$ at the boundary comes precisely
from $\gamma _0(u/d^a)$, the Neumann boundary value (cf.\ Remark 3.11
below).

\subhead 3.2 Decompositions involving systems of traces \endsubhead

 For large $s$, we moreover have representations in terms of
 consecutive sets of traces and Poisson operators. The consecutive
 sets are defined in the following theorem.

\proclaim{Theorem 3.8} Let 
$M$ be a positive integer. 

$1^\circ$ With the Poisson operator  $K_{(0)}$ defined in Theorem {\rm
3.4}, denote
$$
K_{(j)}=\tfrac1{j!}d^jK_{(0)};\tag3.17
$$
it is a Poisson operator of order $-j$ in the Boutet de Monvel calculus.
Then 
$$
\gamma _jK_{(k)}=\delta _{jk}I\text{ for }j\le k,\quad \gamma _jK_{(k)}=\Psi  _{jk}\text{ for }j> k,\tag3.18
$$
where the $\Psi _{jk}$ are $\psi $do's on $\partial\Omega $ of order
$j-k$. With
$$
\varrho _M=\pmatrix \gamma _0\\ \vdots \\ \gamma
_{M-1}\endpmatrix, \quad \Cal K_M=\pmatrix K_{(0)}&\hdots & K_{(M-1)}\endpmatrix,\tag3.19
$$
the composition of $\varrho _M$ and $\Cal K_M$ is a triangular invertible $M\times M$-matrix:
$$
\varrho _M\Cal K_M=\Psi _{(M)}=\pmatrix 1&0&\hdots&0\\ \Psi
_{10}&1&\hdots&0\\
\vdots &\vdots&\ddots&\vdots\\
\Psi _{M-1,0}&\Psi _{M-1,1}&\hdots& 1\endpmatrix,\tag3.20
$$
where $\Psi _{jk}=\delta _{jk}I$ for $j\le k$.
Thus we can define  the right inverse of $\varrho _M$
$$
\widetilde {\Cal K}_M\equiv \Cal K_M\Psi _{(M)}^{-1},\tag3.21
$$
it maps $\prod _{0\le
j<M}B_p^{s-j-1/p}(\partial\Omega )$ into $\overline H_p^{s}(\overline\Omega )$ for
all  $s\in {\Bbb R}$.

$2^\circ$ Define 
$$
K^{\mu }_{(j)}=\tfrac1{\Gamma (\mu +1+j)}d^{\mu }K_{(j)}, \quad j\in{\Bbb N}_0,\tag3.22
$$
and, with $\gamma ^\mu _ju=\Gamma (\mu +1+j)\gamma _j(u/d^\mu )$ (cf. {\rm (1.15)}),
$$
\varrho ^\mu _M=\pmatrix \gamma ^\mu _0\\ \vdots \\ \gamma ^\mu 
_{M-1}\endpmatrix, \quad \Cal K^\mu _M=\pmatrix K^\mu _{(0)}&\hdots & K^\mu _{(M-1)}\endpmatrix,\tag3.23
$$
then $$
\varrho ^\mu _M\Cal K^\mu _M=\Psi _{(M)},\tag3.24 $$
 whereby 
$$
\widetilde{\Cal K}^\mu _M\equiv \Cal K^\mu _M\Psi _{(M)}^{-1} 
\tag3.25 
$$
is a right inverse of $\varrho
^\mu _{M}$.
Here, 
$$
K^\mu _{(j)}\colon B_p^{s-\mu -j-1/p}(\partial\Omega )\to e^+d^{\mu }\overline H_p^{s-\mu }(\Omega )\cap e^+d^{\mu +j}\overline H_p^{s-\mu -j}(\Omega
)
\cap H_p^{(\mu +j)(s)}(\overline\Omega ),
$$ when
$s>\mu +j+\tfrac1p $. Moreover, $\Cal K^\mu _M$ and $\widetilde{\Cal K}^\mu _M$ map 
$\prod _{0\le j<M}B_p^{s-\mu -j-1/p}(\partial\Omega )$
into  \linebreak$ e^+d^\mu \overline H_p^{s-\mu }(\Omega )\cap H_p^{\mu (s)}(\overline\Omega )$ when $s>\mu +M-\tfrac1{p'}$.

\endproclaim

\demo{Proof} $1^\circ$. Since $d$ identifies with $t$ as in Remark 3.1
near $\partial\Omega $, it is verified immediately that $\gamma
_jK_{(k)}=\delta _{jk}$ when $j\le k$. For $j>k$ it is an elementary
fact in the Boutet de Monvel calculus that the composition $\gamma
_jK_{(k)}$ results in a $\psi $do on $\partial\Omega $ of order
$j-k$. These facts lead to (3.20), where the triangular matrix is
clearly invertible (being the sum of the identity matrix and a
nilpotent matrix), as a continuous operator from $\prod _{0\le
j<M}B_p^{s-j-1/p}(\partial\Omega )$ to itself, for all $s$. The 
mapping property of $K_{(j)}$ is well-known for
Poisson operators of order $-j$. 

$2^\circ$. Here the conventions in Remark 3.2 are used; the
traces $\gamma _j^\mu $ are taken from the interior of  $\Omega $, and an extension
$e^+$ is tacitly understood in the definition of $K^\mu _{(j)}$ in
(3.22). Clearly,
$$
\gamma ^\mu _jK^\mu _{(k)}\varphi =\gamma _jK_{(k)}\varphi \text{ for all }j,k.
$$
Then the formula (3.24) follows immediately from (3.20), and (3.25)ff is
a consequence.

That $K^\mu _{(j)}$ maps into $e^+d^{\mu +j}\overline H_p^{s-\mu -j}(\Omega
)\cap H_p^{(\mu +j)(s)}(\overline\Omega )$ follows from Theorem 3.4, since
$K^\mu _{(j)}$ is a constant times $K^{\mu +j}_{(0)}$. It also ranges
in $e^+d^{\mu }\overline H_p^{s-\mu }(\Omega )$, since $d^jK_{(0)}$ is a
Poisson operator of order $-j$. In the
collected statement on $\Cal K^\mu _M$
 and  $\widetilde{\Cal K}^\mu _M$, the space with $j=0$
is common to the mappings.
\qed

\enddemo

Similar calculations hold with $\Omega $ replaced by ${\Bbb R}^n_+$, $K_{(0)}$
replaced by $K_0$. 

\example{Remark 3.9} We here correct some minor flaws in \cite{27},
pages 515--516. The normalizing factor $i^j$ in the definition of
$K_j$ in \cite{27} (1.7) should be removed here, since $D_n^j$ is replaced
by $\partial_n^j$ in the trace definitions in this part of the paper. Moreover, since $\gamma _jK_k$
produces a nontrivial
$\psi $do on ${\Bbb R}^{n-1}$ when $j>k$, the right inverse $\widetilde{\Cal K}_M$
defined above in (3.21) should be used instead of $\Cal K_M$ in all the
formulas. Then formula (5.14) in \cite{27} will contain some more terms
$x_n^{\mu +k}e^+K_0(S_{jk}\gamma ^\mu _{j}u)$ for $0\le k<j$, with
 $\psi $do's $S_{jk}$
 on ${\Bbb R}^{n-1}$ of
order $j-k$. 
These are a minor corrections
that do not change the outcome (5.15) of the theorem. 
\endexample

We can now generalize Theorem 3.4 to sets of traces as follows:

\proclaim{ Theorem 3.10} Let $M\in{\Bbb N}$ and $s>\mu +M-\tfrac1{p'}$. With
 $\varrho ^\mu  _Mu$ and $\widetilde{\Cal K}^\mu _M$ 
defined as in
Theorem {\rm 3.8},
 the elements $u\in H_p^{\mu (s)}(\overline\Omega )$ have unique decompositions 
$$
u=\widetilde{\Cal K}^\mu _M \varphi +w\in e^+d^\mu \overline
H_p^{s-\mu }(\Omega )\cap H_p^{\mu (s)}(\overline\Omega )+  H_p^{(\mu +M)(s)}(\overline\Omega ),\tag3.26
$$
where  $\varphi   $ runs through $\prod _{0\le j<M}B_p^{s-\mu -j-1/p}(\partial\Omega )$ and
$w $ runs through $ H_p^{(\mu +M)(s)}(\overline\Omega )$ (equal to $\dot H_p^s(\overline\Omega )$ if $s-\mu \in\,]M-\tfrac1{p'},M+\tfrac1p [\,$). Here
$\varrho  ^\mu _Mu$ equals  $\varphi   $.

\endproclaim

\demo{Proof} For $u\in H_p^{\mu (s)}(\overline\Omega )$, set $\varphi =\varrho ^\mu 
_Mu$ and $v=\widetilde{\Cal K}^\mu _M \varphi $. Then $w=u-v$ belongs
to $H_p^{\mu (s)}(\overline\Omega )$ and satisfies
$$
\varrho ^\mu _M(u-v)=\varphi -\varphi =0;
$$
hence lies in $H_p^{(\mu +M)(s)}(\overline\Omega )$, in view of (2.12). This
gives a representation of the elements of $H_p^{\mu (s)}(\overline\Omega )$ as
desired. (Extensions by zero $e^+$ are understood.)
An application of $\varrho ^\mu _M$ to (3.26) shows that $\varphi $ must
necessarily  equal $\varrho ^\mu _Mu$, since $\widetilde{\Cal K}^\mu _M$ is a
right inverse of $\varrho ^\mu _M$. 
\qed
\enddemo

It can also be remarked that since
$$
K^\mu _{(j)}=\tfrac1{\Gamma (\mu +1+j)j!}d^{\mu +j}K_{(0)},
$$
we find by setting $\Theta _{\mu ,M}=\bigl(\tfrac1{\Gamma (\mu +1+j)j!}\delta
_{jk}\bigr)_{j,k=0,\dots,M-1}$ and $\psi =\Theta _{\mu ,M}\Phi
^{-1}_M\varphi $ that the term $\widetilde {\Cal K}^\mu _M\varphi$ in
(3.26) satisfies 
$$
\widetilde {\Cal K}^\mu _M\varphi 
={\sum}_{0\le j<M}d^{\mu +j}K_{(0)}\psi _j,\tag3.27
$$
where each $\psi _j$ runs through $B_p^{s-\mu -j-1/p}(\partial\Omega )$.

The progress in this theorem in comparison with Th.\ 5.4 of \cite{27}
is that it gives a precise global parametrization of the elements of the
$\mu $-transmission space $H_p^{\mu (s)}(\overline\Omega )$ in the case of arbitrary
domains $\Omega $, clarifying how the entering elements of the space $e^+d^{\mu }\overline
H_p^{s-\mu }(\Omega )$ look (compare with (2.14)). 
Moreover, it shows explicitly how the
structure of the operators $K^{\mu }_{(j)}$ assures that they map into 
$H_p^{(\mu +j)(s)}(\overline\Omega )$, not just into 
the  spaces 
$e^+d^{\mu +j}\overline H_p^{s-\mu -j}(\Omega
)$ and $ e^+d^{\mu }\overline H_p^{s-\mu }(\Omega
)$ that are less regular over the interior $\Omega $.

Systems of traces have lately been considered  by Abatangelo, Jarohs
and Saldana in \cite{3} in the case of the
fractional Laplacian on the ball, with explicit formulas, and most
recently in \cite{2} (with more coauthors) on the halfspace.

\example {Remark 3.11} When $P$ satisfies Hypothesis 1.1 (recall $a>0$),
then both the spaces $ H_p^{a(s)}(\overline\Omega )$ and   $ H_p^{(a-1)(s)}(\overline\Omega )$ 
are defined, for $s>a-\tfrac1{p'}$ resp.\
$s>a-1-\tfrac1{p'}$. Here $ H_p^{a(s)}(\overline\Omega )$ is for $s>a-\tfrac1{p'}$ the
subspace of elements $u\in  H_p^{(a-1)(s)}(\overline\Omega )$ with $\gamma
^{a-1}_0u=0$, cf.\ (2.12). In fact, $\gamma _0^{a-1}u$ plays the role
of a (nonhomogeneous) Dirichlet boundary {\it value}, and
nonhomogeneous Dirichlet problems for $P$
can be considered on $ H_p^{(a-1)(s)}(\overline\Omega )$ with good solvability
properties,  cf.\ \cite{27}ff. In this context, the
next trace $\gamma _1^{a-1}u$ can be regarded as a Neumann boundary
value; it also enters in nonhomogeneous boundary value problems
(\cite{26,30,32}). It is easy to check that when $u\in  H_p^{a(s)}(\overline\Omega )$, then $$
\gamma _1^{a-1}u=\gamma _0^au.\tag3.28
$$ 
In other words, when $u$ satisfies the homogeneous Dirichlet
condition, the lowest nontrivial boundary value $\gamma _0(u/d^a)$
(or its normalized version $\gamma _0^au$) is the Neumann value.
\endexample

Just like Theorem 3.4 could be generalized to a statement in the
H\"older-Zygmund scale, Theorem 3.10 has such a generalization; it goes as follows:

\proclaim{Theorem 3.12} Let $M\in{\Bbb N}$.
The operators $K^{\mu }_{(j)}$ 
defined  in Theorem {\rm 3.10} map $C_*^{s-\mu -j}(\partial\Omega )$
into $e^+d^{\mu }\overline C_*^{s-\mu }(\Omega )\cap e^+d^{\mu +j}\overline C_*^{s-\mu -j}(\Omega )\cap C_*^{(\mu
+j)(s)}(\overline\Omega )
$ 
when $s>\mu +j$, and  $\widetilde{\Cal K}^\mu _M$ maps 
$\prod _{0\le j<M}C_*^{s-\mu -j}(\partial\Omega )$
into  $  e^+d^{\mu }\overline C_*^{s-\mu }(\Omega )\cap  C_*^{\mu (s)}(\overline\Omega )$ when $s>\mu +M-1$.

 For  $s>\mu +M-1$,  the elements of
 $C^{\mu (s)}_*(\overline\Omega )$ have unique decompositions
$$
u=\widetilde{\Cal K}^\mu _{M}\varphi  +w\in e^+d^\mu \overline
C_*^{s-\mu }(\Omega )\cap C_*^{\mu (s)}(\overline\Omega )+  C_*^{(\mu +M)(s)}(\overline\Omega );
\tag3.29$$
here $w$ 
runs through $
C_*^{(\mu +M)(s)}(\overline\Omega )$ (equal to  $\dot C^s_*(\overline\Omega )$ if $s-\mu \in \,]M-1,M[\,$), and  $\varphi =\varrho ^\mu _Mu$ runs
through $\prod _{j<M}C_*^{s-\mu -j}(\partial\Omega )$.

\endproclaim

Here we are of course primarily interested in the results for noninteger
positive values of the exponents, where the spaces are ordinary
H\"older spaces, $C^s_*=C^s$ for $s\in{\Bbb R}_+\setminus{\Bbb N}$, but the
$C^s_*$ spaces are useful e.g.\ by having good interpolation
properties --- and of course by allowing statements without exceptional
parameters. Moreover, the H\"older-Zygmund spaces allow negative indices;
it is useful to know
(cf.\  \cite{40} and \cite{26}) that the identification of spaces $\dot C^s_*(\overline\Omega )$ and $e^+\overline C_*^s(\Omega
)$ takes place for $-1<s<0$.

The above results will be used with $\mu =a$ in Sections 4 and 5.

\head 4. The regularity of  solutions of fractional Schr\"odinger
Dirichlet problems \endhead

In preparation for the study of heat equation regularity, we shall
consider a related problem for the Schr\"odinger equation, which is of
interest in itself.  
Consider the Dirichlet problem for the fractional Schr\"odinger equation:
$$
r^+Pu+V u=f  \text{ in }\Omega ,\quad \operatorname{supp}u\subset \overline\Omega ,\tag4.1
$$
where the potential  $V$ is a $C^\infty $-function, and $P$ satisfies Hypothesis 1.1.

A regularity discussion of solutions to (4.1) is carried out in
relatively low-order function spaces by Fall \cite{16} for
translation-invariant symmetric operators, under weak smoothness hypotheses on
$P,V$ and $\Omega $. More recently, Diaz, Gomez-Castro and Vazquez
\cite{15} have studied solution properties for quite irregular potentials $V$. 

Recall that for the usual Laplacian $\Delta $, it makes no
difference in the regularity of Dirichlet solutions whether a
$C^\infty $-function $V$ is added or not; the solutions of 
$$
-\Delta u+V u=f  \text{ in }\Omega ,\quad \gamma _0u=0\text{ on }\partial\Omega ,
$$
satisfy for
all $k\in{\Bbb N}_0$, $ 0<\delta  <1$:
$$
f\in \overline C^{k+\delta }(\Omega ) \iff u\in  \overline C^{2+k+\delta }(\Omega
).
$$
(Cf.\ e.g.\ Courant and Hilbert \cite{14} p.\ 349; $V$ just enters as a zero-order term.)

In contrast, for noninteger powers of $-\Delta $, and operators $P$
satisfying Hypothesis 1.1 with $a \notin{\Bbb N}$, the regularity may be considerably restricted in
comparison with the case $V =0$. This is linked to the fact that the
multiplication by a nonzero function $V $ does not fit into the symbol
sequence $p\sim
\sum_{j\in{\Bbb N}_0}p_j$ with $p_j(x,-t\xi )=t^{2a-j}(-1)^jp_j(x,\xi )$.

We first improve the regularity as far as we can by using the known
regularity results for the Dirichlet problem
$$
r^+Pu=g \text{ in }\Omega ,\quad \operatorname{supp}u\subset \overline\Omega .\tag4.2
$$

 \proclaim{ Theorem 4.1} Let $P$ satisfy Hypothesis {\rm 1.1} for
 some  $a\in {\Bbb R}_+\setminus{\Bbb N}$, and let $V\in \overline C^\infty (\Omega )$.  Let $u\in \dot
 H^a (\overline\Omega )$ satisfy {\rm (4.1)} for some $f\in \overline C^\infty (\Omega )$. Then 
$u\in C^{a(3a)}_*(\overline\Omega )$.  

The conclusion also holds if merely $f\in \overline C^a (\Omega )$. In fact,
the solutions $u$ with  $f\in \overline C^a (\Omega )$ run through $
C^{a(3a)}_*(\overline\Omega )$, with $\gamma _0^a u$ running through
$C_*^{2a}(\partial\Omega )$.

In cases where $f\in \overline C^s(\Omega )$ for some $s\in [0,a[\,$,  $u\in C_*^{a(2a+s)}(\overline\Omega )$.

\endproclaim

\demo{Proof} 
By a variant of the variational treatment of
$r^+P$ (adding the term $(Vu,v)$ to $Q_0(u,v)$), the operator $r^+P+V$
is Fredholm 
from $\{u\in \dot H^a (\overline\Omega )\mid r^+Pu+Vu\in L_2(\Omega )\}=\{u\in \dot H^a (\overline\Omega )\mid r^+Pu\in L_2(\Omega )\}=D(P_{\operatorname{Dir},2})
$ to
$L_2(\Omega )$.

Let $u\in \dot H^a (\overline\Omega )$ be a solution of
(4.1) with $f\in \overline C^\infty (\Omega )$.
Using the regularity theory for (4.2), 
we shall improve the knowledge of the regularity of $u$ in a finite
number of iterative steps, as in a related situation in \cite{28},
pf.\ of Th.\ 2.3:

Recall
the well-known general embedding properties for $p,p_1\in \,]1,\infty [\,$:
$$
\dot H_p^{a}(\overline\Omega )\subset e^+L_{p_1}(\Omega ),\text{ when }\tfrac
n{p_1}\ge \tfrac n p - a,  \quad \dot H_p^{a}(\overline\Omega )\subset
\dot C^0(\overline\Omega  )\text{ when }a>\tfrac np.\tag4.3 
$$

We make a finite number of iterative steps to reach
the information $u \in \overline C^0(\Omega )$, as follows: 
Define $p_0,p_1, p_2,\dots$, with $p_0=2$ and $q_j=\frac n {p_j}$ for
all the relevant $j$,
such that
$$
q_j=q_{j-1}-a \text{ for }j\ge 1.
$$
 This means that $q_j=q_0-ja$; we stop the sequence at $j_0$ the first
 time we reach a  $q_{j_0}\le 0$. 

As a first step, we note that $u \in \dot
 H^a(\overline\Omega )\subset e^+L_{p_1}(\Omega )$ implies that $f-V u\in
 L_{p_1}(\Omega )$, whence $u \in
 H_{p_1}^{a(2a)}(\overline\Omega )$ by
 \cite{27} Th.\ 4.4 applied to $r^+Pu\in  L_{p_1}(\Omega )$. Then
  $u\in \dot
 H^a_{p_1}(\overline\Omega )$ in view of (2.11). 
In the next step we use the embedding $\dot
 H^a_{p_1}(\overline\Omega )\subset e^+L_{p_2}(\Omega )$ to conclude in a
 similar way that
 $u \in \dot H^a_{p_2}(\overline\Omega )$, and so on. If $q_{j_0}<0$, we have
 that $\frac n{p_{j_0}}<a$, so $u \in \dot
 H^a_{p_{j_0}}(\overline\Omega )\subset \dot C^0(\overline\Omega )$. If $q_{j_0}=0$, the
 corresponding $p_{j_0}$ would be $+\infty $, and we see at least that
 $u \in e^+L_{p}(\Omega )$ for any large $p$; then one step more
 gives that $u \in \dot C^0(\overline\Omega )$.
 
The rest of the argumentation relies on H\"older estimates. 
By the
regularity results recalled in Section 2, Theorem 2.3ff.,
$$
f-V u \in \overline C^0(\Omega) \implies u \in C^{a(2a)}_*(\overline\Omega )\subset
e^+d^a \overline C^{a}(\overline\Omega )+\dot C^{2a-\varepsilon  }(\overline\Omega )\subset
\dot C^{a}(\overline\Omega ).\tag4.4
$$ 
Next, $f-V u \in \overline C^{a}(\Omega)$ implies $
u \in C^{a(3a)}_*(\overline\Omega )$.       

Clearly, only the smoothness
$f\in \overline C^{a}(\Omega)$ is needed for the whole argumentation.

Conversely, if $u\in C^{a(3a)}_*(\overline\Omega )$, then  $r^+Pu\in \overline C_*^{2a}(\Omega
)\subset \overline C^a(\Omega )$, and since $u\in \dot C_*^{3a-\varepsilon
}(\overline\Omega )+e^+d^a\overline C_*^{2a}(\Omega )\subset e^+\overline C^a(\Omega )$,
also $Vu\in \overline C^a(\Omega )$. Here $\gamma _0^au$ runs through
$C_*^{2a}(\partial\Omega )$ (cf.\ Theorem 3.6).

For the last assertion in the theorem, let $f\in \overline C^s(\Omega
)$. Since $\overline C^s(\Omega )\subset \overline C^0(\Omega )$, the inductive
steps in the above proof can be carried through until we reach (4.4).
Now since $s<a$, $f-Vu\in \overline C^s(\Omega )$, so by Theorem 2.3,
$u\in   C^{a(2a+s)}_*(\overline\Omega )$.
\qed

\enddemo

\example{Remark 4.2} If the potential $V$ vanishes near
 $\partial\Omega $,
the regularity of $u$ when $f\in \overline C^\infty (\Omega )$ can be lifted all
the way to $\Cal E_a(\overline\Omega )$, by relying on the interior regularity
$H^{a(s)}(\overline\Omega ) \subset H^{s}_{\operatorname{loc}}(\Omega )$: Let
$V\in C_0^\infty (\Omega )$, $f\in \overline C^\infty (\Omega )$.
 Then $u\in H^{a(s)}(\overline\Omega )$ implies
$Vu\in \dot H^{s}(\overline\Omega  )$, hence $r^+Pu=f-Vu\in \overline H^s(\Omega )$
and consequently $u\in H^{a(s+2a)}(\overline\Omega )$. So the regularity can be lifted
in steps of size $2a$, eventually reaching $u\in
\bigcap_sH^{a(s)}(\overline\Omega )=\Cal E_a(\overline\Omega )$.
(One can instead argue that the operator $P+V$ in this case satisfies the
$a$-transmission condition at $\partial\Omega $, since $V$ does not
enter in boundary patches close to $\partial\Omega $, so the results
from \cite{27} hold for $P+V$.) 
The regularity
connected with this kind of
potentials was used in Ghosh, Salo and Uhlmann \cite{21}.

\endexample

 We shall now see that the general regularity results stop at this
 level. Namely,  when $V$ is
nonvanishing on part of the boundary, a higher regularity than in
Theorem 4.1 can only
hold if $\gamma ^a_0u$ vanishes on that part. 

We show this for $0<a<1$;
cases of higher $a\notin{\Bbb N}$ can be treated in similar ways.

\proclaim{Theorem 4.3} Let $P$ be as in Theorem {\rm 4.1} with
$0<a<1$, and let $V\in  \overline C^\infty (\Omega )$. 
Let $u\in C_*^{a(3a)}(\overline\Omega )$ and $f\in \overline C^{a }(\Omega
)$ solve {\rm (4.1)}. 

$1^\circ$ Asssume that  $1/V\in  \overline C^\infty (\Omega )$.
If  $u\in C_*^{a(3a+\delta )}(\overline\Omega )$ 
and $f\in \overline C^{a+\delta }(\Omega )$ 
for some $\delta >0$, then
$\gamma ^a_0u=0$.

$2^\circ$ Let $V(x)\ne 0$ for $x$ in an open subset $\Sigma $ of the boundary $\partial\Omega $.
If  $u\in C_*^{a(3a+\delta )}(\overline\Omega )$ 
and $f\in \overline C^{a+\delta }(\Omega )$ 
for some $\delta >0$, then
$\gamma ^a_0u=0$ on $\Sigma $.

\endproclaim

\demo{Proof}
$1^\circ$. First consider the case  where $0<a<\frac12$. Then
$2a<1$, and we can take $\delta \in \,]0,1-2a[\,$, so that also $2a+\delta <1$. Assume
that  $u\in C_*^{a(3a+\delta )}(\overline\Omega )$ and $f\in \overline C^{a+\delta
}(\Omega )$. 
By (3.16) with  $s=3a+\delta $, $\mu =a$,
$$
u=K^a_{(0)}\gamma ^a_0u +w_0,\quad w_0\in C_*^{(a+1)(3a+\delta
)}(\overline\Omega )=\dot C_*^{3a+\delta }(\overline\Omega ),
$$
where we used that $3a+\delta -a\in\,]0,1[\,$, cf.\ the last statement
in Theorem 3.6. Here $K^a_{(0)}\gamma
^a_0u=d^az$ with $z=\frac1{\Gamma (a+1)}e^+K_{(0)}\gamma ^a_0u\in e^+\overline C^{2a+\delta }(\Omega )$. Since $\dot
C_*^{3a+\delta }(\overline\Omega )\subset d^a\dot C_*^{2a+\delta }(\overline\Omega )$,
we have that
$$
u/d^a=z+w'_0, \quad w_0'=w_0/d^a\in \dot C_*^{2a+\delta }(\overline\Omega ) , 
$$
and hence, in local coordinates at the boundary, where $d$ is replaced by $x_n$,
$$
\underline u(x',x_n)=x_n^a(\underline z(x',0)+O(x_n^{2a+\delta
}))=x_n^a\underline z(x',0)+O(x_n^{3a+\delta })\text{ for small }x_n>0.
$$
On the other hand, $u=V ^{-1}(f-r^+Pu)\in \overline C^{a+\delta
}(\Omega )$ and therefore has an expansion in local coordinates
$$
\underline u(x',x_n)=\underline u(x',0)+O(x_n^{a+\delta })\text{ for
small }x_n>0.
\tag4.5
$$
Comparing the two expansions, we first conclude that $\underline
u(x',0)=0$, and next, that $\underline z(x',0)=0$.
This shows that $\gamma ^a_0u$ must be 0. 

Now consider the case $a\in [\frac12,1[\,$. Here
$2a\in [1,2[\,$, and we consider $\delta >0$ satisfying $\delta <1-a$;
then also $\delta
<2-2a$, so that $2a+\delta \in \,]1,2[\,$ and $a+\delta <1$. Assume
that  $u\in C_*^{a(3a+\delta )}(\overline\Omega )$ and $f\in \overline C^{a+\delta
}(\Omega )$. 
By (3.30) with $M=2$, $s=3a+\delta $,
$$
u=K^a_{(0)}\gamma ^a_0u +K^{a}_{(1)}(\gamma ^a_1u -\Psi _{10}\gamma ^a_0u)+w_1,\text{ with }w_1\in C_*^{(a+2)(3a+\delta
)}(\overline\Omega )=\dot C_*^{3a+\delta }(\overline\Omega ),
$$
where we used that $3a+\delta -a-2\in \,]-1,0[\,$. Here $$
K^a_{(0)}\gamma ^a_0u +K^{a}_{(1)}(\gamma ^a_1u -\Psi _{10}\gamma
^a_0u)=d^az+d^{a+1}z_1,
$$
with $z=\frac1{\Gamma (a+1)}e^+K_{(0)}\gamma ^a_0u\in e^+\overline
C^{2a+\delta }(\Omega )$ and $z_1=\frac1{\Gamma
(a+1)}e^+K_{(1)}(\gamma ^a_1u-\Psi _{10}\gamma ^a_0u)\in   e^+d\overline
C^{2a+\delta -1}(\Omega )$. Then
$$
u/d^a=z+dz_1+w'_1, \quad w_1'=w_1/d^a\in \dot C_*^{2a+\delta }(\overline\Omega ) , 
$$
and hence, in local coordinates where $d$ is replaced by $x_n$,
$$
\underline u(x',x_n)=x_n^a\underline
z(x',0)+O(x_n^{a+1})
+O(x_n^{3a+\delta})
\text{ for small }x_n>0.\tag4.6
$$
By comparison with the expansion (4.5)
we can again first conclude that $\underline
u(x',0)=0$, and next that $\underline z(x',0)=0$, so we find again that
$\gamma ^a_0u$ must be 0. This shows $1^\circ$.

For $2^\circ$, we just carry the above argumentation through in
coordinate patches intersecting the boundary in open subsets $\Sigma '$ of
$\Sigma $ with $\overline {\Sigma '}\subset \Sigma $.  
\qed

\enddemo

As a corollary we find for the resolvent equation:

\proclaim{Corollary 4.4} Let $P$ satisfy Hypothesis {\rm 1.1} with
$0<a<1$, and let $\lambda \ne 0$. 
Let $u\in C_*^{a(3a)}(\overline\Omega )$ and $f\in \overline C^{a }(\Omega
)$ solve 
$$
r^+Pu-\lambda  u=f  \text{ in }\Omega ,\quad \operatorname{supp}u\subset \overline\Omega .\tag4.7
$$ 
If  $u\in C_*^{a(3a+\delta )}(\overline\Omega )$ 
and $f\in \overline C^{a+\delta }(\Omega )$ 
for some $\delta >0$, then
$\gamma ^a_0u=0$.
\endproclaim

\demo{Proof} This is the special case of Theorem 4.3 $1^\circ$ where
$V=-\lambda $. \qed
\enddemo

Expressed in words: When $f\in \overline C^{a+\delta }(\Omega )$, we
cannot have a solution $u$ of (4.7) in the corresponding solution space $C_*^{a(3a+\delta )}(\overline\Omega )$ for (4.2)
without imposing an extra boundary condition $\gamma ^a_0u=0$.

 Recall from
Remark 3.11 that for $u\in C_*^{a(t)}(\overline\Omega )$, $\gamma _0^au$ can be
regarded as the Neumann boundary value. 

The vanishing of $\gamma _0^au$ occurs in \cite{15} in a different
context, namely when $V$ is very irregular at $\partial\Omega $; such
solutions are called $a$-flat there.

\example{Remark 4.5}
What is it that happens when $s$ in the parameter $2a+s$ passes from
$a$ to $a+\delta $, $\delta >0$? Recall that we are dealing with the operator
$r^+P$ going from $E^s_1=C_*^{a(2a+s)}(\overline\Omega  )$ to $E^s_0= \overline C^s_*(\Omega
)$, say when $s> 0$.
Here we always have that  (with $\varepsilon $ active
if $2a+s\in{\Bbb N} $)
$$
E^s_1\subset \dot C_*^{2a+s(-\varepsilon )}(\overline\Omega )+e^+d^a\overline
C_*^{a+s}(\Omega )\subset d^a[\dot C^a(\overline\Omega )+e^+\overline C^a(\Omega
)]\subset \dot C^a(\overline\Omega ).\tag4.8 
$$
But note also that (cf.\ (2.19))
$$
E^s_1=C_*^{a(2a+s)}(\overline\Omega  )\supset \Cal E_a(\overline\Omega )=e^+d^a\overline C^\infty
(\Omega ),
$$ and for $\delta >0$ there are elements of $e^+d^a\overline C^\infty
(\Omega )$ not lying in $e^+\overline C^{a+\delta }(\Omega )$ (for example
$d^a$ itself).
So for $s=a+\delta $, $E^s_1$ contains nontrivial elements of $\dot
C^a(\overline\Omega  )\setminus e^+ \overline C^{a+\delta }(\Omega )$. Briefly expressed, 
$$
\aligned
E^s_1&\subset E^s_0 \text{ when }s\le a,
\\
E^s_1&\not\subset  E^s_0 \text{ when }s> a.
\endaligned\tag4.9
$$
The inclusion $E^s_1\subset E^s_0$ is necessary in order to define a
resolvent acting in $E^s_0$; this is not possible for $s>a$.

\endexample

We shall also see what happens for larger $s$, in particular when $s$
grows to $\infty $. By $[s]$ we denote the largest integer $\le s$.

\proclaim{Theorem 4.6} Let $P$ and $V$ be as in Theorem {\rm 4.3},
and let  $V(x)\ne 0$ for $x$ in an open subset $\Sigma $ of
$\partial\Omega $.  Let $s_0>a$, and let $k$ be the largest integer
$<s_0-a$.

Assume that $u$ and $f$ solve {\rm (4.1)} with $u\in
C_*^{a(2a+s_0)}(\overline\Omega  )$, $f\in \overline C_*^{s_0}(\Omega )$. Then 
$$
\gamma ^a_ju=0\text{ on $\Sigma $ for }j=0,1,\dots, k.\tag4.10
$$

In particular, if $V=-\lambda $ for a $\lambda \in {\Bbb C}\setminus
\{0\}$, {\rm (4.10)} holds on $\partial\Omega $.
 
\endproclaim

\demo{Proof} 
Let $s_0-a-k=\delta _0$, it is $>0$. The assumptions on $u$ and $f$
hold also with $s_0$ replaced by $t=a+k+\delta $, any $\delta \in
\,]0,\delta _0]$, where we can choose $\delta $ as small as we please.

{\it The case $0<a<\frac12$.} Take $\delta $ so small that $2a+\delta
<1$, a fortiori $a+\delta <1$; then $a+t=2a+k+\delta $ is noninteger, with $[a+t]=[2a+k+\delta
]=k$. Moreover, $t$ is noninteger with  $[t]=[a+k+\delta ]=k$. We shall apply Theorem 3.12 to 
 $ C_*^{a(2a+t)}(\overline\Omega  )$, with
$\mu ,s$ in the theorem  replaced by $a,2a+t$. Note that when we
define $M=[a+t+1]$, then $a+t<M<a+t+1$, i.e., $M-1<a+t<M$. Moreover, $M-1=k$. Here $a+t$
corresponds to $s-\mu $ in Theorem 3.12. The theorem then applies to
show that $u\in C_*^{a(2a+t)}(\overline\Omega  )$ satisfies
$$
u\in e^+d^a\overline C^{a+t}(\Omega)+ \dot C^{2a+t}(\overline\Omega ). \tag4.11
$$
In local coordinates at the boundary (as described in Remark
3.1), where $d$ is replaced by $x_n$, this implies that with $\underline v=\underline u/x_n^a$,
$$
\align
\underline u (x',x_n)&=x_n^a[\underline v (x',0)+x_n\partial_n\underline v (x',0)+\dots\tag4.12\\
&\qquad +\tfrac
1{k!}x_n^{k}\partial_n^{k}\underline v (x',0)+O(x_n^{a+t})]+O(x_n^{2a+t})\\
&=x_n^a\underline v (x',0)+x_n^{a+1}\partial_n\underline v (x',0)+\dots+\tfrac
1{k!}x_n^{a+k}\partial_n^{k}\underline v (x',0)+O(x_n^{2a+t}),
\endalign 
$$
for $x_n\to 0+$, by Taylor expansion of $\underline v $. Here the functions $\partial_n^{j}\underline v (x',0)$ are
proportional to the traces $\gamma ^a_j\underline u $. 

On the other hand, when $x\in\Sigma $, there is a neighborhood $U$ of
$x$ in ${\Bbb R}^n$ such that $V^{-1}$ exists as a $C^\infty $-function on
$U\cap\overline\Omega $, and then $u=V^{-1}(f-r^+Pu)$ is $C^t$ on
$U\cap\overline\Omega $. Then in local coordinates, since $[t]=k$,
$$
\underline u (x',x_n)=\underline u (x',0)+x_n\partial_n\underline u (x',0)+\dots+\tfrac
1{k!}x_n^{k}\partial_n^{k}\underline u (x',0)+O(x_n^{a+k+\delta })\tag4.13
$$
there.
Note here that $0<a<1<a+1<2<\dots <k-1<a+k-1<k<a+k$.
Comparing the two expansions (4.12) and (4.13), we find successively
that $\underline u (x',0)=0$, hence $\underline v (x',0)=0$, hence $\partial_n\underline u (x',0)=0$,
hence $\partial_n\underline v (x',0)=0$, etc., until we reach the information 
$\partial_n^{k}\underline u (x',0)=0$. Then the term $\tfrac
1{k!}x_n^{a+k}\partial_n^{k}\underline v (x',0)$ must also vanish, since what
is left in (4.13) is $O(x_n^{a+k+\delta })$. This shows that the traces $\gamma ^a_ju $
must vanish on $U\cap \partial\Omega $ for $j\le k$.

{\it The case $\frac12\le a<1$.} Take $\delta $ so small that $2a+\delta
\in \,]1,2[\,$ and $a+\delta <1$, then $a+t=2a+k+\delta $ is noninteger, with $[a+t]=[2a+k+\delta
]=k+1$. Moreover, $t$ is noninteger with  $[t]=[a+k+\delta ]=k$. Here,
when we apply Theorem 3.12 to 
 $ C_*^{a(2a+t)}(\overline\Omega  )$, the number $M=[a+t+1]$ equals $k+2$. Then
 the inclusion (4.11) implies, in local coordinates,
$$
\underline u (x',x_n)=x_n^a\underline v (x',0)+x_n^{a+1}\partial_n\underline v (x',0)+\dots+\tfrac
1{(k+1)!}x_n^{a+k+1}\partial_n^{k+1}\underline v (x',0)+O(x_n^{2a+t}).\tag4.14
$$
 Again we compare this with (4.13); this leads immediately to the 
vanishing of the
 terms with factors up to $x_n^{a+k-1}$ and $x_n^k$, and then the
 estimate $O(x_n^{a+k+\delta })$ in (4.13) shows that also
 $\partial_n^{k}\underline v (x',0)=0$.
 We conclude that the traces $\gamma ^a_ju $ vanish on $U\cap \partial\Omega $ for $j\le k$.
\qed
\enddemo

There is also a result in the $C^\infty $-category:

\proclaim{Theorem 4.7} Let $P$ and $V$ be as in Theorem {\rm 4.3},
with $1/V\in \overline C^\infty (\Omega )$. If $u\in \Cal E_a(\overline\Omega )$ and
$f\in \overline C^\infty (\Omega )$ satisfy {\rm (4.1)}, then $u\in \dot
C^\infty (\overline\Omega )$ (all traces vanish).
\endproclaim

\demo{Proof} This is a corollary to Theorem 4.6, when we let
$s\to\infty $ there. But it is easier to see directly: In view of
(1.4), $u$ must lie in
$\Cal E_a(\overline\Omega )\cap e^+\overline C^{\infty }(\Omega)$, which can only hold when all
boundary values vanish (as also noted in \cite{28}).\qed
\enddemo

 We leave it to the reader
to formulate the last result when $V$ only vanishes on part of the boundary.

\head 5. The regularity of  solutions of fractional heat
Dirichlet problems \endhead

 The results on the resolvent equation can now be applied in a discussion of the 
regularity of solutions of the fractional heat equation.

\proclaim{Theorem 5.1}  Let $P$ satisfy Hypothesis {\rm 1.1} with
$0<a<1$. When $u\in  $ \linebreak $\overline
W^{1,1}({\Bbb R};C_*^{a(3a)}(\overline\Omega))$, it satifies
$$
\aligned
r^+Pu(x,t)+\partial_tu(x,t)&=f(x,t)\text{ on }\Omega \times {\Bbb
R},\\
u(x,t)&=0 \text{ for }x\notin\overline\Omega  .
\endaligned
\tag5.1
$$
with $f(x,t)\in
L_1({\Bbb R};\overline C^{a }(\Omega))$; here   $\gamma
^a_0u\in  \overline W^{1,1}({\Bbb R};C_*^{2a}(\partial\Omega))$ can take any value.

However, if for some $\delta >0$,
 $u(x,t)\in  \overline
W^{1,1}({\Bbb R};C_*^{a(3a+\delta )}(\overline\Omega))$ and $f(x,t)\in
L_1({\Bbb R};\overline C^{a +\delta }(\Omega))$, then   $\gamma
^a_0u=0$.

\endproclaim

\demo{Proof}
For the first statement, we note that $\partial_tu\in L_1({\Bbb
R};C_*^{a(3a)}(\overline\Omega)) \subset $ \linebreak $ L_1({\Bbb
R};\dot C^{a}(\overline\Omega))$, so $r^+Pu+\partial_tu\in  L_1({\Bbb
R};\overline C^{a}(\Omega)) $ as asserted.

 For the second statement, we need only consider a small $\delta >0$ with $a+\delta <1$.
The functions are sufficiently regular to allow
 Fourier transformation with respect to $t$, leading to the equation
(where we denote $\Cal F_{t\to\tau }g(x,t)=\grave g(x,\tau )$):
$$
r^+P\grave u(x,\tau )+i\tau \grave u(x,\tau )=\grave f(x,\tau )\text{
on }\Omega \times {\Bbb R},\quad \grave u(x,\tau )=0\text{ for
}x\notin \overline\Omega  .\tag5.2
$$
By assumption, 
$$
u\in \overline W^{1,1}({\Bbb R};C_*^{a(3a+\delta )}(\overline\Omega )),
\;
f\in 
L_1({\Bbb R};\overline C^{a+\delta }(\Omega ))  ,
\;
\gamma _0^au\in  \overline W^{1,1}({\Bbb
R};C_*^{2a+\delta }(\partial\Omega)),
\tag5.3
$$ 
so since the Fourier transform maps $L_1({\Bbb R};X)$ into $C^0({\Bbb
R};X)$, the ingredients in  (5.2) are in spaces: 
$$
 r^+P\grave u\in C^1({\Bbb
R};\overline C^{a+\delta }(\Omega)),\;\tau \grave u(x,\tau )\in
C^0({\Bbb R};C_*^{a(3a+\delta )}(\overline\Omega )),\; \grave f\in C^0({\Bbb
R};\overline C^{a+\delta }(\Omega)),\tag5.4
$$ 
and the equation holds pointwise in $\tau $ (this use of Fourier
transformation of functions valued in Banach spaces is justified by
the analysis in Amann \cite{4}).
 
At each $\tau \ne 0$ we can apply Corollary 4.4 to (5.2), to see that
if $\delta >0$, then $\grave u(x,\tau )$
cannot be in $C_*^{a(3a+\delta )}$ in $x$ unless $\gamma ^a_0\grave
u=0$ for that value of $\tau $.

Observing this at all $\tau \ne 0$, we see in view of the continuity
in $\tau $ that if $\delta >0$, $\gamma ^a_0\grave
u(x,\tau )$ vanishes. By the injectiveness of the Fourier transform, also $\gamma ^a_0
u(x,t)$ vanishes. \qed
\enddemo

Again, expressed in words:  When $f$ is $\overline C^{a+\delta }(\Omega )$ in $x$, $u$
cannot be in the corresponding solution space $C_*^{a(3a+\delta
  )}(\overline\Omega )$ with respect to $x$
without satisfying the extra boundary condition $\gamma ^a_0u=0$.

As a corollary, we have a similar result for $t$ in a finite interval:

\proclaim{Corollary 5.2} With $I=\,]0,T[\,$, let $u\in  \overline
W^{1,1}(I;C_*^{a(3a)}(\overline\Omega))$ with $u(x,0)=0$; then it satifies
$$
\aligned
r^+Pu(x,t)+\partial_tu(x,t)&=f(x,t)\text{ on }\Omega \times I,\\
u(x,t)&=0 \text{ for }x\notin\overline\Omega  ,\\
u(x,0)&=0,
\endaligned
\tag5.5
$$
with $f(x,t)\in
L_1(I;\overline C^{a }(\Omega))$ (here  $\gamma
^a_0u\in  \overline W^{1,1}({I};C_*^{2a}(\partial\Omega))$ can be
freely prescribed for positive $t$).

However, if for some $\delta >0$,
 $u(x,t)\in  \overline
W^{1,1}(I;C_*^{a(3a+\delta )}(\overline\Omega))$ and $f(x,t)\in $
 \linebreak $
L_1(I;\overline C^{a +\delta }(\Omega))$, then   $\gamma
^a_0u=0$.

\endproclaim

\demo{Proof} First extend $u$ and $f$ by 0 for $t<0$, and next extend
the resulting functions across $t=T$ by reflection in $t$. This
results in  functions $\tilde u$
resp.\ $\tilde f$ that satisfy the hypotheses of Theorem 5.1, so the
conclusions from that theorem carry over. \qed
\enddemo

For higher regularity classes, one finds that an increasing number of
traces must vanish in order to have a solution of (5.5). The following
theorem is proved analogously to Theorem 5.1 and Corollary 5.2, using Theorem 4.6:

\proclaim{Theorem 5.3} Let $P$ satisfy Hypothesis {\rm 1.1} with
$0<a<1$,  let $s>a$,  and let $k$ be the largest integer $<s-a$. For $I=\,]0,T[\,$, let 
 $u(x,t)\in  \overline
W^{1,1}(I;C_*^{a(2a+s )}(\overline\Omega))$ and $f(x,t)\in
L_1(I;\overline C_*^{s}(\Omega))$ solve the problem {\rm (5.5)}. Then   $\gamma
^a_ju=0$ for all $j\le k$.

\endproclaim

We also have a conclusion in the $C^\infty $-category, found by
applying the information in Theorem 5.3 for all $s$ (using that
$\Cal E_a(\overline\Omega)=\bigcap_sC_*^{a(2a+s )}(\overline\Omega)$ and 
$\overline C^\infty (\Omega)= \bigcap_s\overline
C_*^s(\Omega)$):

\proclaim{Corollary 5.4} Let $P$ satisfy Hypothesis {\rm 1.1} with
$0<a<1$. Assume that $u(x,t)\in  \overline
W^{1,1}(I;\Cal E_a(\overline\Omega))\equiv \bigcap_s\overline
W^{1,1}(I;C_*^{a(2a+s )}(\overline\Omega))
$ and $f(x,t)\in
L_1(I;\overline C^\infty (\Omega))\equiv \bigcap_sL_1(I;\overline
C_*^s(\Omega))$ solve the problem {\rm (5.5)}. Then $\gamma ^a_ju$
vanishes for all $j\in{\Bbb N}_0$.

\endproclaim

It should be noted that the regularity theorem of Ros-Oton and Vivas
\cite{50}
shows that for solutions of (5.5),
$$
f\text{ is }C^a \text{ in $x$ and $C^{\frac12 }$ in $t$}
\implies u/d^a
\text{ is }C^{2a }\text{ in $x$ and $C^{1}$ in
$t$},\tag5.6
$$
if $a\ne \frac12$ (with a slightly weaker statement for $a=\frac12$,
cf.\ (1.7));
this is consistent with the first assertion in Corollary 5.2. But an
extension of the upper indices from
$a$ to $a+\delta $, resp.\  $2a$ to $2a+\delta $, would possibly need restrictive hypotheses as in the
second assertion.

\medskip

Received XXXX 2018; revised XXXX 2018.
\medskip

 {\it E-mail address}: grubb\@math.ku.dk 


\Refs



\ref\key 1  \by N.\ Abatangelo \paper\rm Large s-harmonic functions and
boundary blow-up solutions for the fractional Laplacian \jour\it
 Discrete Contin.\ Dyn.\ Syst.\ \vol 35 \yr2015\pages 5555--5607
\endref 


\ref\key 2 \by N. Abatangelo, S. Dipierro, M. M. Fall, S. Jarohs, A. Saldana 
\finalinfo arXiv:1806.05128 
\paper\rm
Positive powers of the Laplacian in the half-space under Dirichlet boundary conditions
\endref
 
\ref\key 3  \by N. Abatangelo, S. Jarohs and A. Saldana \paper\rm 
Integral representation of solutions to higher-order
fractional Dirichlet problems on balls \jour\it
  Comm. Contemp. Math. \finalinfo to appear, arXiv: 1707.03603 
\endref 


\ref\key 4  \by H. Amann \paper\rm  Operator-valued Fourier multipliers, vector-valued Besov spaces, and
 applications\jour\it Math. Nachr. \vol 186 \yr 1997 \pages 5--56\endref



\ref\key 5 \by U. Biccari, M. Warma and E. Zuazua \paper\rm Local
regularity for fractional heat equations \finalinfo arXiv: 1704.07562
\endref


\ref\key 6  \by B. M. Blumenthal and R. K. Getoor \paper\rm The
asymptotic distribution of the eigenvalues for a class of Markov
operators \jour\it Pacific J. Math. \vol 9\yr 1959 \pages 399--408
\endref
 
\ref\key 7 \paper\rm     Censored stable processes
 \by   K. Bogdan, K. Burdzy and Z.-Q. Chen \vol 127\pages 89--152
\jour\it Prob. Theory Related Fields\yr 2003 
\endref

\ref\key 8 \by    M. Bonforte, Y. Sire and J. L. Vazquez 
\paper\rm     Existence, uniqueness and asymptotic behaviour for
fractional porous medium equations on bounded domains \jour\it
Discrete Contin. Dyn. Syst. \vol 35 \yr 2015 \pages 5725-–5767 
 \endref




\ref\key 9  \by L.\ Boutet de Monvel\paper\rm Boundary problems for pseudo-differential
operators \jour\it Acta Math.\ \vol 126 \yr1971 \pages 11--51\endref
 
\ref\key 10  \by L.\ Caffarelli and  L.\  Silvestre\paper\rm An extension problem related
to the fractional Laplacian \jour\it Comm.\ Part.\ Diff.\ Eq.\ \vol 32 \yr
2007 \pages 1245--1260
\endref



\ref\key 11 
\by H. Chang-Lara and G. Davila \paper\rm Regularity for solutions of non
local parabolic equations\jour\it Calc. Var. Part. Diff.
Equations \vol 49 \yr2014 \pages 139-–172
\endref

\ref\key 12 \by Z.-Q. Chen and R. Song \paper\rm Estimates on Green
functions and Poisson kernels for symmetric stable processes \jour\it
Math. Ann. \vol 312 \yr1998 \pages 465--501
\endref 



\ref\key 13  \by R. Cont and P. Tankov \book Financial Modelling with
Jump Processes, 
Financial Mathematics Series \publ
Chapman \& Hall/CRC \publaddr Boca Raton, FL \yr 2004
\endref


\ref\key 14  \by R. Courant and D. Hilbert \book Methods of
Mathematical Physics II \publ Interscience Publishers \publaddr New
York \yr 1962
\endref



\ref\key 15  \by J.I.Diaz, D. Gomez-Castro and J.L. Vazquez \paper\rm
The fractional Schr\"odinger equation with general nonnegative
potentials. The weighted space approach \finalinfo arXiv:1804.08398
\endref


\ref\key 16  \by M.M. Fall \paper\rm
Regularity estimates for nonlocal Schr\"odinger equations \finalinfo arXiv:1711.02206
\endref

\ref\key 17 
\by M. Felsinger and M. Kassmann \paper\rm Local regularity for parabolic
nonlocal operators\jour\it Comm. Part. Diff. Equations \vol 38
\yr2013 
\pages 1539-–1573\endref

\ref\key 18 \by M. Felsinger, M. Kassmann and P. Voigt \paper\rm The
Dirichlet problem for nonlocal operators  \jour\it Math. Z.\vol 279\yr
2015 \pages779-–809
\endref

\ref\key 19  
\paper\rm Regularity theory for general stable operators: parabolic equations
\by
X. Fernandez-Real and X. Ros-Oton\jour\it
J. Funct. Anal. \vol 272 \yr2017 \pages 4165--4221
\endref


\ref\key 20 \by R. Frank and L. Geisinger \paper\rm Refined
semiclassical asymptotics for fractional powers of the Laplace
operator\jour\it J. Reine Angew. Math. \vol712 \yr 2016 \pages 1-–37
\endref





\ref\key 21  \by T. Ghosh, M. Salo and G. Uhlmann \paper\rm The
Calder\'on problem for the fractional Schr\"odinger equation
\finalinfo arXiv:1609.09248
\endref

\ref \key 22 \by M. Gonzalez, R. Mazzeo and Y. Sire \paper\rm Singular
solutions of fractional order conformal Laplacians \jour\it
J. Geom. Anal.\vol 22 \yr2012\pages 845-–863 \endref

\ref\key 23  \by G.\ Grubb\paper\rm Pseudo-differential boundary problems in Lp spaces,
Comm.\ Part.\ Diff.\ Eq.\ \vol 15 \yr1990 \pages 289--340\endref 


 \ref\key 24 \by 
{G.~Grubb}\book Functional Calculus of Pseudodifferential
     Boundary Problems.
 Pro\-gress in Math.\ vol.\ 65, Second Edition \publ  Birkh\"auser
\publaddr  Boston \yr 1996
\endref

\ref\key 25 \by G. Grubb\book Distributions and Operators. Graduate
Texts in Mathematics, 252 \publ Springer \publaddr New York\yr 2009
 \endref


\ref\key 26  \by G.\ Grubb\paper\rm  
Local and nonlocal boundary conditions for $\mu $-transmission
and fractional elliptic pseudodifferential operators \jour\it 
 Analysis and P.D.E.\  \vol 7 \yr 2014\pages 1649--1682\endref

\ref\key 27  \by G.\ Grubb\paper\rm Fractional Laplacians on domains, 
a development of H\"o{}rmander's theory of $\mu$-transmission
pseudodifferential operators \jour\it Adv.\ Math.\  \vol 268 \yr 2015 \pages
478--528\endref

\ref\key 28  \by G.\ Grubb\paper\rm Spectral results for mixed problems
and fractional elliptic operators \jour\it J. Math. Anal. Appl.  \vol 421 \yr 2015 \pages
1616--1634\endref

\ref\key 29  \by  G.\ Grubb \paper\rm Regularity of spectral
fractional Dirichlet and Neumann problems \jour\it Math.\
Nachr.\  \vol 289 \yr 2016\pages 831--844
\endref


\ref\key 30  \by  G.\ Grubb \paper\rm Integration by parts and  Pohozaev
identities for space-dependent fractional-order operators \jour\it J.\
Diff.\ Eq.\ \vol 261 \yr 2016\pages 1835--1879
\endref


\ref\key 31 \by G. Grubb \paper\rm  Regularity in $L_p$ Sobolev spaces
of solutions to fractional heat equations \jour\it J. Funct. Anal.\vol
274
\yr 2018 \pages 2634--2660
\endref

\ref\key 32 \by G. Grubb \paper\rm   Green's formula and a Dirichlet-to-Neumann operator for
fractional-order pseudodifferential operators \jour\it
Comm. Part. Diff. Equ.
\yr 2018 \finalinfo 
doi.org/10.1080/03605302.2018.1475487, arXiv:1611.03024 
\endref


\ref\key 33 \by G. Grubb \paper\rm  Fractional-order operators: boundary
problems, heat equations  \inbook in \it
Mathematical Analysis and Applications
--- Plenary Lectures, ISAAC 2017, Vaxjo Sweden, Springer Proceedings
in Mathematics and Statistics \eds L. G. Rodino and
J. Toft\publ Springer\publaddr Switzerland \pages 51--81 \yr 2018 
\endref



\ref\key 34  \by W. Hoh and N. Jacob \paper\rm On the Dirichlet problem
for pseudodifferential operators generating Feller semigroups
\jour\it J. Functional Anal. \vol 137 \yr 1996 \pages 19--48 \endref

\ref\key 35  \by L.\ H\"o{}rmander\paper\rm Seminar notes on
pseudo-differential operators and boundary problems 
\finalinfo Lectures at IAS Princeton 1965-66, 
available from Lund University 
https://lup.lub.lu.se/search/
\endref


\ref\key 36  \by L.\ H\"ormander\book The Analysis of Linear Partial
Differential Operators, III \publ Springer Verlag \yr 1985 \publaddr
Berlin
\endref

\ref\key 37 \by
N. Jacob \book Pseudo Differential Operators and Markov
Processes. Vol. I--3
\publ Imperial College Press \publaddr London \yr 2001 \endref

\ref\key 38  \by T. Jakubowski \paper\rm The estimates for the Green function in Lipschitz
 domains for the symmetric stable processes \jour\it
 Probab. Math. Statist. \vol 22
 \yr 2002 \pages  419-441\endref

\ref\key 39 \by 
T. Jin and J. Xiong \paper\rm Schauder estimates for solutions of linear
parabolic integro-differential equations \jour\it Discrete
Contin. Dyn. Syst. \vol 35 \yr2015 \pages 5977-–5998\endref

\ref\key 40  \by J. Johnsen \yr 1996 \jour\it Math. Scand.\paper\rm Elliptic
boundary problems and the Boutet de Monvel calculus in Besov and
Triebel-Lizorkin spaces \vol 79  \pages 25--85\endref


\ref\key 41  \by T. Kulczycki \paper\rm Properties of Green function of symmetric stable
 processes \jour\it Probab. Math. Statist. \vol 17 \yr1997 \pages
 339-364 \endref

\ref\key 42 \by N. S. Landkof\book Foundations of Modern Potential
Theory.
Die
Grundlehren der mathematischen Wissenschaften, Band 180 \bookinfo  (Translated from the Russian by A. P. Doohovskoy.) \publ
Springer-Verlag
\publaddr New York-Heidelberg \yr 1972 \endref



\ref\key 43 \by T. Leonori, I. Peral, A. Primo and F. Soria \paper\rm
Basic estimates for solutions of a class of nonlocal elliptic and
parabolic equations
\jour\it
Discrete Contin. Dyn. Syst. \vol 35 \yr 2015 \pages 6031--6068
\endref



\ref\key 44  \by R. Musina and A.I. Nazarov \paper\rm On fractional
Laplacians \jour\it Comm. Part. Diff. Eq.\vol 39 \yr2014 \pages
1780-–1790
\endref

\ref\key 45 \paper\rm
Nonlocal elliptic equations in bounded domains: a survey
\by X. Ros-Oton
\jour\it Publ. Mat. \vol60 \yr2016 \pages 3--26\endref

\ref\key 46 \paper\rm Boundary  regularity,  Pohozaev identities and
nonexistence results
\by X. Ros-Oton \inbook Recent Developments in Nonlocal Theory \publ
De Gruyter \publaddr Berlin
\pages 335--358
\endref

\ref\key 47  \by X.\ Ros-Oton and J.\ Serra\paper\rm The Dirichlet problem for the
fractional Laplacian: regularity up to the boundary \jour\it J.\ Math.\ Pures
Appl.\ \vol 101 \yr 2014 \pages  275--302\endref


\ref\key 48  \by  X.\ Ros-Oton and J.\ Serra \paper\rm
The Pohozaev identity for the fractional Laplacian \jour\it
Arch.\ Rat.\ Mech.\ Anal.\ \vol 213 \yr 2014 \pages 587--628\endref



\ref\key 49  \by  X.\ Ros-Oton and J.\ Serra \paper\rm
Regularity theory for general stable operators \jour\it J. Differential Equations \vol 260 \yr
2016 \pages 8675-–8715\endref



\ref\key 50  \by X.\ Ros-Oton and H.\ Vivas \paper\rm Higher-order
boundary regularity estimates for nonlocal parabolic equations
\jour\it Calc. Var. Partial Differential Equations, Art. 111,
\vol  57   \yr 2018 \page 20 pp \endref



\ref\key 51 \by R. Servadei and E. Valdinoci
\paper\rm On the spectrum of two different fractional operators
\jour\it Proc. Roy. Soc. Edinburgh \vol 144 \yr 2014 \pages 831--855
\endref

\ref\key 52 \by L. Silvestre \paper\rm Regularity of the obstacle problem
for a fractional power of the Laplace operator \jour\it Comm. Pure
Appl. Math. \vol 60 \yr 2007 \pages 67--112
\endref

\ref\key 53 \by M.~E. Taylor\book
 Pseudodifferential Operators \publ
Princeton University Press \publaddr Princeton, NJ \yr1981
\endref



\endRefs
\enddocument